\documentclass[12pt]{article}
\usepackage{amsfonts,amsmath}
\usepackage[not]{luecke}
\newtheorem{Definition}{Definition}
\newtheorem{Proposition}{Proposition}
\newtheorem{Theorem}{Theorem}
\newtheorem{Corollary}{Corollary}
\newtheorem{Lemma}{Lemma}
\newenvironment{Proof}{\hfill\par\noindent{\bf Proof:}}{\hfill\rule{2mm}{2mm}\vspace{0.1cm}\par\noindent}
\newenvironment{Remark}{\hfill\par\noindent{\bf Remark.}}{\hfill\par\noindent}
\newenvironment{Example}{\hfill\par\noindent{\bf Example.}}{\hfill\par\noindent}

\begin{document}

\def\eps{\varepsilon}
\def\d{\mathrm{d}}
\def\DD{\mathcal{D}}
\def\AA{\mathcal{A}}
\def\BB{\mathcal{B}}
\def\QQ{\mathcal{Q}}
\def\SS{\mathcal{S}}
\def\HH{\mathcal{H}}
\def\M{\mathbb{M}}
\def\Z{\mathbb{Z}}
\def\N{\mathbb{N}}
\def\R{\mathbb{R}}
\def\CC{\mathcal{C}}
\def\MM{\mathcal{M}}
\def\PP{\mathcal{P}}
\def\LL{\mathcal{L}}
\def\SS{\mathcal{S}}
\def\TT{\mathcal{T}}
\def\ii{\lbrack 0,1 \rbrack}
\def\SBV{\mathrm{SBV}}
\def\S{\mathrm{S}}
\def\P{\mathrm{P}}
\def\T{\mathrm{T}}
\def\E{\mathrm{E}}
\def\D{\mathrm{D}}
\def\C{\mathrm{C}}
\def\MM{\mathcal{M}}
\def\UU{\mathcal{U}}
\def\VV{\mathcal{V}}
\def\MS{\mathrm{MS}}
\def\BZ{\mathrm{BZ}}
\def\BV{\mathrm{BV}}
\def\ms{\mathrm{ms}}
\def\bz{\mathrm{bz}}
\def\sing{\mathrm{sing}}
\def\supp{\mathrm{supp}}
\def\dom{\mathrm{dom}}
\def\disc{\mathrm{disc}}
\def\argmin{\mathrm{argmin}}
\def\kn{{\scriptstyle \frac{\kappa}{n}}}
\def\kpn{{\scriptstyle \frac{\kappa+1}{n}}}
\def\kpnzwei{{\scriptstyle \frac{\kappa+2}{n}}}
\def\kmn{{\scriptstyle \frac{\kappa-1}{n}}}
\def\koeff{{\scriptstyle \frac{n}{\mu^2}}}
\def\koenn{{\scriptstyle \frac{n}{\mu^2_n}}}
\def\koess{{\scriptstyle \frac{n}{\mu^2_s}}}
\def\RI{\lbrack {\scriptstyle \frac{n-1}{n}},1\rbrack}

\title{Segmentation of Time Series: Parameter Dependence of Blake-Zisserman and Mumford-Shah
Functionals and the Transition from Discrete to Continuous}

\author{A. Kempe$^1$\thanks{Partially supported by DFG
 Graduate Programme
`Applied Algorithmic Mathematics' at the TU M\"{u}nchen and DFG grant SFB 386 at the LMU M\"{u}nchen}\and V. Liebscher$^1$\and G. Winkler$^1$ \and
O. Wittich$^2$
\\ \\
{\small $^1$Institute of Biomathematics and Biometry}
\\ {\small GSF-National Research
Center for Environment and Health} \\ {\small Ingolst\"{a}dter Landstr.1, D - 85764 Neuherberg} \\ {\texttt{\small
kempe,liebscher,winkler@gsf.de}}\\
\\ {\small $^2$M 12, Centre for Mathematical Science, TU M\"{u}nchen}\\
{\small Boltzmannstr. 3, D - 85747 Garching bei M\"{u}nchen}
\\ {\texttt{\small wittich@ma.tum.de}}}
\date{}
\maketitle
\newpage
\begin{abstract}The paper deals  with variational approaches to the segmentation of time series into smooth pieces,
but allowing for sharp breaks. In discrete time, the corresponding functionals are of Blake-Zisserman type.  Their natural counterpart in
continuous time are the Mumford-Shah functionals. Time series which minimise these functionals are proper estimates or representations of
the signals behind recorded data. We focus on consistent behaviour of the functionals and the estimates, as parameters vary or as the
sampling rate increases.

For each time continuous time series  $f\in L^2 (\lbrack 0,1\rbrack)$ we take conditional expectations w.r.t. to
$\sigma$-algebras generated by finer and finer partitions of the time domain into intervals, and thereby construct a sequence
$(f_n)_{n\in\N}$ of discrete time series.  As $n$ increases this amounts to sampling the continuous time series with more and more
accuracy.

Our main result is consistent behaviour of segmentations w.r.t. to variation of parameters and increasing sampling rate.
\end{abstract}

\vskip0.5cm


\vskip0.5cm


\vskip0.5cm
\newpage

\tableofcontents
\newpage

\section{Introduction}

\luecke[Muss Generelles aus Abstract hier nochmal kommen? Ich denke JA.]

We will first introduce and discuss the adopted concept of segmentation, and  define the functionals we are dealing with. Then the main
result will be stated, and finally, we sketch the plan of the paper.

\subsection{Segmentations}

A fundamental {\em variational ansatz} for the  segmentation of time series is the minimisation of functionals which penalise undesired
properties of the estimate against fidelity to the data. The latter will be measured by the $L^2$-distance of the estimate and data. The
penalty should depend on
\begin{enumerate}
\item the measure of the set of `jumps' or `breaks' between regions of `smoothness',
\item a notion of
smoothness which restricts the behaviour of the signal between two subsequent jumps.
\end{enumerate}
The notion of `jump' and `smoothness' will be made precise shortly. Let us illustrate the concepts by way of two examples. Suppose that the
{\em segmentation}
is a square integrable function
$g$ on the closed unit interval $U =\lbrack 0,1\rbrack$. As segmentations, we
allow functions $f = t + F$ with a {\em right continuous step function} $t$ and a Sobolev differentiable function
$F$. These functions are of {\em special bounded variation}, see Section \ref{SBV}. The {\em Mumford-Shah functionals}
in time dimension one are defined as
\begin{equation*}
\MS_{\gamma,\mu,g}(f) := \gamma j(t) +
    \frac{1}{\mu^2}\int_U \vert f^{\prime}\vert^2 ds + \Vert g -
    f\Vert^2
\end{equation*}
where $j(t)$ is the {\em number of discontinuities} of $t$ in
$(0,1)$,  $f^{\prime}$ is the Sobolev derivative of $F$, and the parameters $\gamma$ and $\mu$ control the number of breaks and the degree
of smoothness. This is the one-dimensional version of the functionals introduced in \cite{MumShah:85} and \cite{MumShah:89} for space
dimension two.

In this model, {\em jumps} are the discontinuities of $t$ which by {\em Sobolev's embedding theorem} can be identified with the
discontinuities of
$f$. `Smoothness' is measured by the $L^2$-norm
of the derivative of $F$ within the intervals between subsequent jumps; it coincides there with the derivative of
$f$. {\em Fidelity to the data}, finally, is measured by the
$L^2$-distance of  the segmentation to data.

The second example -- which is in fact a special case of the Mumford - Shah functional if $\mu = 0$ -- is called {\em Potts - functional},
inspired by the functional introduced in \cite{Potts:52} as a generalization of the Ising model (\cite{Ising:25}) in statistical mechanics.
It is defined for step functions $t$ and given by
\begin{equation*}
\P_{\gamma,g}(t) := \gamma j(t) + \Vert t - g\Vert^2 .
\end{equation*}
The notions of jump and of fidelity to the data is the same as for the Mumford - Shah functional, whereas the notion of smoothness is
considerably stronger, only functions which are even constant on the intervals between consecutive jumps are possible outcomes of the
segmentation procedure. The Mumford - Shah functional (and hence the Potts-functional as well) can be extended to $L^2 (U)$ in the sense
that it is given by the expression above, if the equivalence class of $f\in L^2 (U)$ contains a function which can be written as $f = t +
F$, and that it is equal to $\infty$, if not (cf. Definition \ref{DefMS}).

In contrast to the continuous setting, there are no `obvious' notions of smoothness or jumps for discrete time - series
$\underline{g} = (g^0,...,g^{n-1})\in\R^n$. One possibility is to
consider those points $\kappa\in \lbrace 0,...,n-2\rbrace$ as jumps, where the difference $\vert f^{\kappa +1} - f^{\kappa}\vert$ exceeds a
given threshold (cf. Definition
\ref{discdef}). One particular functional, where this notion is
immanent is the {\em Blake - Zisserman} functional (see
\cite{Blake:83}, \cite{BlZi:87}). For $\gamma,\mu \geq 0$ it is
given by
\begin{equation*}
    \BZ_{\gamma,\mu,\underline{g}} (\underline{f}) := \sum_{\kappa =0}^{n-2}
    \min\lbrace \vert f^{\kappa +1} - f^{\kappa}\vert^2 /\mu^2,\gamma\rbrace + \sum_{\kappa =0}^{n-1} (f^{\kappa} -
    g^{\kappa})^2.
\end{equation*}
An equivalent definition (cf. Lemma \ref{edgemod}) focussing more closely on the nature of jumps (´edges') in this functional by showing
that the ´min' in the penalizing term arises from a minimization with respect to a larger space of covariables including explicit
information about ´jumps', was given in
\cite{GemanGeman:84}. As for the Mumford - Shah functional, the
special case $\mu = 0$, i.e.
\begin{equation*}
    \P_{\gamma,\underline{g}} (\underline{f}) :=
    \gamma \,\vert \lbrace \kappa = 0,...,n-1\,:\, f^{\kappa +1}\neq f^{\kappa}\rbrace\vert
    + \sum_{\kappa =0}^{n-1} (f^{\kappa} -
    g^{\kappa})^2
\end{equation*}
is called {\em (discrete) Potts model}. It also favors segmentations $\underline{f}$ which are constant between consecutive jump points.

In the variational ansatz, we consider as segmentations those functions or vectors, respectively, which minimize the given functionals.
That immediately leads to the following questions:
\begin{itemize}
\item[{\rm (i)}] Is there always a minimizer of the given
functional ? \item[{\rm (ii)}] How does it depend on the (free) parameters $\gamma$ and $\mu$ ?
\end{itemize}
In view of the discussion about jumps and smoothness above, we may also ask a question about the relation between the discrete and the
continuous setup:
\begin{itemize}
\item[{\rm (iii)}] Is there an embedding of the discrete situation
into the continuous one, i.e. a way of discrete sampling from a continuous `truth', such that the discrete segmentations converge to the
segmentations of the continuous signal ?
\end{itemize}
One exact formulation of question {\rm (iii)} and an affirmative
answer about the relation of Mumford - Shah and Blake - Zisserman
functional is the main subject of this paper and is provided by
Theorem \ref{Hauptsatz} in the next subsection. As a byproduct of
our analysis, we obtain in that case answers to questions (i) and
(ii) as well.

\subsection{The Main Result}\label{1.2}

For an exact formulation of the third question, we have to specify
the embedding of the discrete into the continuous situation. Since
we are considering signals $g\in L^2 (U)$, it is not suitable to
consider vectors arising from the evaluation of a signal $g$ at
several distinct time points. Instead, we consider its {\em
conditional expectation} with respect to the $\sigma$- algebra
generated by a fixed {\em partition} of $U=\lbrack 0,1\rbrack$.
The intuition behind this procedure is, that the output of our
measuring device is truly an average of the signal over a short
period of time. Sampling with more and more accuracy thus means to
decrease the length of these periods. We adopt the following
conventions:

\begin{Definition}\label{SetUp} Let $n\in\N$ and the {\em equidistant setup} be
given by
\begin{equation*}
    \SS (n) := \lbrace \kn \,: \kappa = 0,...,n\rbrace\subset U
\end{equation*}
and $\sigma_n$ denote the $\sigma$-algebra
\begin{equation*}
\sigma_n := \sigma \left(\lbrack 0,{\scriptstyle
\frac{1}{n}}),\lbrack {\scriptstyle \frac{1}{n}},{\scriptstyle
\frac{2}{n}}),..., \lbrack {\scriptstyle
\frac{n-2}{n}},{\scriptstyle \frac{n-1}{n}}),\lbrack {\scriptstyle
\frac{n-1}{n}},1\rbrack\right).
\end{equation*}
The conditional expectation of a signal $g\in L^2 (U)$ with
respect to $\sigma_n$ will be shortly denoted by $g_n :=
E(g\,\vert\,\sigma_n)$.
\end{Definition}

\begin{Remark} The fact that the rightmost interval is closed has
no further significance. We found that it is most convenient to
deal with the right boundary point in the way that we consider
only step functions which are continuous at $s=1$.
\end{Remark}

As an appropriate discrete input for the Blake - Zisserman
functional we consider now the signal $\underline{g}_n= \pi_ng
:=(g^0_n,...,g^{n-1}_n)\in\R^n$ given by
\begin{equation*}
    g_n^{\kappa} := E (g\,\vert\,\sigma_n) (\kn) =
    n\int_{\kn}^{\kpn} g(s) ds .
\end{equation*}
A one - sided inverse to this {\em discretization map} is provided
for $\underline{f}\in\R^n$ by the step function
\begin{equation*}
    \tau_n{\underline{f}}(s) := f^{n-1}\chi_{\RI}(s) + \sum_{\kappa =0}^{n-2} f^{\kappa} \,\chi_{\lbrack\kn,\kpn )}(s).
\end{equation*}

\begin{Remark} Note that for the conditional expectation of a
signal $g\in L^2 (U)$, we have $g_n = E(g\,\vert\,\sigma_n) =
\tau_n\circ\pi_n g$. This fact will be used frequently in the
sequel.
\end{Remark}

As a consequence of our treatment of the boundary point $s=1$ --
which implies that all $\sigma_n$-measurable functions are right-
continuous step functions with left limits which are additionally
continuous at $s=1$ -- we will assume this property as well for
all step functions considered in this paper. Hence we arrive at
the following definition.

\begin{Definition}\label{Treppen} (1) We denote by $\TT (U)$ the set of all
right-continuous step functions on $U$ with left limits which are
additionally continuous at $1\in U$. (2) We denote by $\TT_n
(U)\subset \TT (U)$ the set of $\sigma_n$-measurable functions on
$U$.
\end{Definition}

With these conventions in mind, we can now define an embedded
version of the Blake - Zisserman functional and of the discrete
Potts - Functional on $L^2 (U)$. Let $\underline{f}_n := \tau_n f
= (f^0,..., f^{n-1})\in\R^n$ given by $f^{\kappa} := f(\kn)$ if
$f\in \TT_n(U)$ and $\underline{g}_n$ given by the conditional
expectation above. Then
\begin{equation*}
    \BZ^n_{\gamma,\mu,g} (f) :=
    \left\lbrace\begin{array}{ll}
\BZ_{\gamma,\mu/n,\underline{g}_n/\sqrt{n}} (\pi_n f/\sqrt{n})&
\mathrm{if} \,
f\in\TT_n(U)  \\
\infty & \mathrm{else}
    \end{array}\right.
\end{equation*}
and
\begin{equation*}
\P^n_{\gamma,g} (f) := \left\lbrace\begin{array}{ll}
\P_{\gamma,\underline{g}_n/\sqrt{n}} (\pi_nf/\sqrt{n}) &
\mathrm{if} \,
f\in\TT_n(U)  \\
\infty & \mathrm{else}
    \end{array}\right.,
\end{equation*}
respectively (cf. Definition \ref{BZD}).

\noindent To be more precise, we embed the functionals defined
above into a family of functionals depending on two real and one
additional rational parameter which reflects the transition from
discrete to continuous. Let
\begin{equation*}
    \T:= \lbrace 1/n\,:\,n\in\N\rbrace\cup\lbrace 0 \rbrace \subset\R
\end{equation*}
be equipped with its relative topology as a subset of $\R$, i.e.
$0$ is the only accumulation point. Consider the (pseudo-) cube
$\QQ = \R^+_0\times \R^+_0 \times \T$. To each point of the cube
corresponds one functional as follows:

\noindent Let $(\gamma, \mu, t)\in\QQ$ and $g\in L^2(U)$. We adopt
the convention that $t=1/n$, $0=1/\infty$, respectively. Then, the
functions $F (\gamma,\mu,t):L^2 (U)\to\R$ are given by
\begin{equation*}
F (\gamma,\mu,t) := \left\lbrace
\begin{array}{ll}
\BZ_{\gamma,\mu,g}^{n} & \gamma > 0, \mu > 0, t > 0 \\
\MS_{\gamma,\mu,g} & \gamma > 0, \mu > 0, t = 0 \\
\P_{\gamma,g}^n & \gamma > 0, \mu = 0, t > 0 \\
\P_{\gamma,g} & \gamma > 0, \mu = 0, t = 0 \\
\d_g^n & \gamma = 0, \mu \geq 0, t > 0 \\
\d_g & \gamma = 0, \mu \geq 0, t = 0 \\
\end{array}\right..
\end{equation*}
Here $\BZ_{\gamma,\mu,g}^n$, $\MS_{\gamma,\mu,g}$ are the {\em
Mumford- Shah-} and {\em Blake-Zisserman functionals} which were
shortly discussed above and which precise form extended to
functionals on $L^2 (U)$ is given in the Definitions \ref{DefMS},
\ref{BZD}, below. The {\em Potts-functionals $\P_{\gamma,g}^n $
and $P_{\gamma,g}$} are as above, the {\em discrete distance
functional} $\d_g^n: L^2 (U)\to\R$ is given by
\begin{equation*}
    \d_g^n (f) :=
    \left\lbrace\begin{array}{ll}\frac{1}{n}\sum_{\kappa=0}^{n-1}(f-g_n)^2(\kn)& f\in\TT_n(U) \\
    \infty & \mathrm{else}\end{array}\right.
\end{equation*}
and finally the {\em continuous distance functional} $\d_g : L^2
(U)\to\R$ given by
\begin{equation*}
    d_g (f) := \int_U (f-g)^2 ds.
\end{equation*}
With these definitions, the main result of this paper reads as
follows:

\begin{Theorem}\label{Hauptsatz}
For the family $F(q):L^2 (U)\to \R,q\in\QQ$, the following
statements hold:
\begin{enumerate}
\item For all $q\in \QQ$ there is a minimizer of $F(q)$. \item Let
$q_s = (\gamma_s,\mu_s,t_s)$, $s\in\N$ converge to
$q=(\gamma,\mu,t)$ in $\QQ$ as $s$ tends to infinity. Denote by
$f^*_{s,s\in\N}$ a sequence of minimizers of $F(q_s)$. Then:
\begin{enumerate}
\item Every convergent subsequence of $f^*_{s,s\in\N}$ converges
to a minimizer of $F(q)$. \item $f^*_{s,s\in\N}$ always contains a
convergent subsequence.
\end{enumerate}
\end{enumerate}
\end{Theorem}

That means, segmentation of a time series using the variational
ansatz with these functionals behaves consistently under variation
of the parameters and under sampling the true signal with more and
more accuracy.

\subsection{Plan of the Paper}

The proof of our main result is mainly based on the fact, that the
minimisation of the functionals under consideration can be
splitted up: Once a set of jumps is fixed, there is a minimizer of
the functional among all admissible functions whose jump set is
contained in the given one. This minimizer can be computed rather
explicitly. The minimization of the functional over all possible
segmentations reduces thus to the minimization with respect to all
jump sets, which will be identified with {\em partitions} of $U$
or, respectively, its associated minimizers. In the sequel, this
common feature of the functionals will be called the {\em
reduction principle}. The minimizers associated to a fixed jump
set will be called {\em partition solvers}.

In the following section, we begin by fixing the notion of {\em
partition} which is basically the exact manifestation of ´jump
set', together with some of its properties as finite sets.
Furthermore, we make precise our point of view by defining what we
understand as a {\em segmentation}. As its most important
manifestation in this paper, we consider {\em functions of special
bounded variation}. In Section \ref{RP}, we give an exact
statement of the {\em reduction principle} and compute the {\em
partition solvers}. Continuity properties of the latter ones are
proved in the subsequent section. The first important implication
of these considerations is relative compactness of the set of
minimizers derived in Section \ref{CP}. As explained at the
beginning of Section \ref{GC}, $\Gamma$-convergence together with
this compactness result implies Theorem \ref{Hauptsatz}. To
establish the $\Gamma$-convergence result Theorem \ref{GammaCube}
in the final section, the basic ingredient is an associated result
of independent interest, given by the Lemmas \ref{liminf} and
\ref{limsup}, on the interchangeability of Sobolev differentiation
and approximation by step functions. This is proved in the
remainder of Section \ref{GC}. We end up with a (then) short proof
of our main result Theorem \ref{Hauptsatz}.

\section{Segmentations and Partitions}

In this section, we introduce what we mean by the segmentation of
a signal which depends on one parameter, for instance a discrete
or continuous time-series.

\subsection{The Continuous Case}

In that case, the signal depends on a continuous set of
parameters. We restrict ourselves to square integrable signals.

\begin{Definition}\label{segclass} Let $U := \ii$ and $L^2 (U)$ denote the Hilbert
space of (equivalence classes) of square integrable functions on
$U$ with respect to Lebesgue measure. A {\em signal} is a function
$g\in L^2 (U)$. A {\em segmentation class} on $U$ is a class $\S
(U)$ of (equivalence classes) of right-continuous functions with
left limits, i.e. there is an injective map $\S (U)\subset\D (U)$
where $\D (U)$ denotes the {\em Skorohod space} (see e.g.
\cite{Billingsley:1999}).
\end{Definition}

\noindent We have the following examples for segmentation classes:

\begin{Example}(1) $\S(U)=\D(U)$. (2) $\S(U)=\TT (U)$, the space of
all right-continuous step functions with left limits which are even continuous at $1\in U$. (3) $\S(U)=\SBV_2 (U)$, the
functions of {\em 2-bounded special variation} discussed below.
\end{Example}

Intuitively, a {\em segmentation} provides a decomposition of a
signal into homogeneous parts which are separated by abrupt
changes ({\em jumps}). The decomposition idea is represented by
the concept of {\em partitions}. For reasons of technical
convenience, we decided that all partitions contain the boundary
points of $I$.

\begin{Definition} Let $\vert M \vert$ denote the cardinality of a subset $M\subset U$.
(1) The {\em partitions} of the interval $U$ are given by the set
\begin{equation*}
    \PP (U) := \lbrace p\subset\ii\,:\, \vert p \vert <\infty ;
    \,0,1\in p\rbrace .
\end{equation*}
(2) The {\em special closed subsets} of the interval $I$ are given
by the set
\begin{equation*}
    \CC (U) := \lbrace c\subset\ii\,:\, c \,\,\mathrm{closed};
    \,0,1\in c\rbrace .
\end{equation*}
(3) The partition $p(f)$ {\em associated to} a given $f\in\D(U)$
with finite set of discontinuities $\disc (f)$ is $p(f) := \disc
(f)\cup \lbrace 0,1\rbrace$.
\end{Definition}
In the sequel, the points in $\disc (f)$ are frequently called
{\em jumps}.

\subsection{The Discrete Case}\label{2.2}

In the discrete case, we still think of a continuous signal in the
sense defined above. In contrast to the continuous case, we can
only sample its values at finitely many (equidistant) time points.
Hence, partitions into continuous parts separated by
discontinuities make no longer sense and we have to substitute it
by something else. Recall the definitions \ref{SetUp} and
\ref{Treppen} of the equidistant setup and the step functions
associated to it.

\begin{Definition}\label{setup} By
\begin{equation*}
    \PP_n (U) :=\lbrace p\in\PP (U)\,:\, p\subset \SS (n)\rbrace
\end{equation*}
we denote the partitions compatible with the equidistant setup
$\SS (n)$.
\end{Definition}

So far, we did not say anything about a decomposition of the
signal into more and less homogeneous parts. In the continuous
case, this was achieved with the help of a partition associated to
the segmentation function. We will do the same now by introducing
a suitable threshold.

\begin{Definition}\label{discdef} A {\em threshold} is a function $T: \SS (n)\to
\R^+_0$. The {\em $T$-partition} of a function $f\in\TT_n(U)$ is
given by
\begin{equation*}
    p_{T}(f):= \lbrace\kn : \vert f(\kpn)-f(\kn)\vert > T (\kappa) ; \kappa =1,...,n-1\rbrace
    \cup \lbrace 0,1\rbrace \in \PP_n (U).
\end{equation*}
Points in $p_{T}(f)$ are called {\em discrete discontinuities} of
$f$ (with respect to $T$).
\end{Definition}

In that sense, we consider those points as points of discontinuity
where the difference of values at adjacent sampling points exceeds
a given threshold. Note that the threshold may as well be adapted,
i.e. depend on the function under consideration.

\begin{Example}(1) For $T = 0$, all points where $f(\kpn)\neq f(\kn)$ are discrete
discontinuities. This choice corresponds to the consideration of
the minimizers of the Potts functional, which are constant off the
jump set. (2) For the Blake-Zisserman functional to parameters
$\gamma,\mu,n$ (see below), the appropriate choice is $T(\kappa)
=\mu\,\sqrt{\gamma/n}$.
\end{Example}

\subsection{Partitions and Hausdorff Metric}

A topology on the set of closed subsets of $U$ is provided by the
Hausdorff metric.

\begin{Definition} Let $c,c^{\prime}$ be closed non-empty subsets of $U$.
Then the distance
\begin{equation*}
    d_H (c,c^{\prime}) := \max \left\lbrace \max_{x\in c} \min_{y\in
    c^{\prime}} \vert x - y \vert , \max_{x\in c^{\prime}} \min_{y\in
    c} \vert x - y \vert \right\rbrace
\end{equation*}
is called {\em Hausdorff distance}.
\end{Definition}

The properties of the Hausdorff distance are summarized by the
following proposition, the proof follows from standard facts,
available for instance in \cite{A:Bee:93} or \cite{Matheron:1975}.

\begin{Proposition}\label{toppartitions}(1) The closed subsets of $U$, provided with
Hausdorff distance, form a compact (hence complete) metric space. (2) The set $\CC (U)$ is a closed subspace, hence as
well compact. (3) The subset $\PP (U)\subset \CC (U)$ is dense, i.e. $\overline{\PP (U)}= \CC (U)$. (4) The subset
$\PP_n (U)\subset \CC (U)$ is finite, hence as well closed and compact.
\end{Proposition}
In the sequel, we will need another characterization of Hausdorff convergence of partitions focussing on the
decomposition of $U$ into intervals.

\begin{Definition} Let $p\in \PP(U)$. By $p^{\mathrm c}:=U-p$ we denote the {\em
complement} of $p$ in $U$. The complement of $p = \lbrace 0=x_0<
x_1<...<x_m =1\rbrace$ is a disjoint union of open subintervals
$\Theta\subset U$. The collection of these subintervals will be
denoted by
\begin{equation*}
    \iota (p) := \lbrace \Theta\subset p^{\mathrm c}\,:\, \Theta
    = ( x_k, x_{k+1})\,:\, k= 0,...,m-1\rbrace .
\end{equation*}
\end{Definition}
Now the characterization of Hausdorff convergence in terms of the
intervals reads as follows.

\begin{Lemma}\label{convchar} Let $p_n\in\PP(U)$ converge to
$p=\lbrace 0=x_0< x_1<...<x_m =1\rbrace\in\PP(U)$ in Hausdorff
metric. Then we have
\begin{enumerate}
\item For all $\Theta\in \iota(p)$ there is a sequence
$\Theta_n\in\iota (p_n)$ such that $\overline{\Theta_n}$ converges
to $\overline{\Theta}$ in Hausdorff metric. \item
$\overline{\Theta_n}= \lbrack a_n,b_n\rbrack$ converges to
$\overline{\Theta}=\lbrack a,b \rbrack$ if and only if $a_n\to a$
and $b_n\to b$. \item Let $\Theta,\Theta^{\prime}\in\iota (p)$ two
adjacent intervals, i.e $\Theta =(x,y)$, $\Theta^{\prime} =
(y,z)$. Let $\Theta_n = (x_n,y_n)$, $\Theta^{\prime}_n =
(y_n^{\prime},z_n^{\prime})$ the sequences of intervals from (i).
Then $\lbrack y_n,y_n^{\prime}\rbrack\to\lbrace y \rbrace$ in
Hausdorff metric and $ y_n,y_n^{\prime}\to y$ uniformly for all
adjacent $\Theta,\Theta^{\prime}\in\iota(p)$.
\end{enumerate}
\end{Lemma}

\begin{Proof} Let $\delta_0 := 1/3\,\min_{x\neq y\in p} \vert x - y\vert$ and $\Theta = (a,b)\in \iota (p)$ be fixed. Without
loss of generality, assume that $n$ is so large that $d_H (p_n,p)
< \delta_0$. (ii) By $a < b$, $a_n < b_n$ and by the assumption,
we have necessarily $\vert a - a_n \vert, \vert b - b_n\vert <
\delta_0$ and thus by definition of the Hausdorff metric
\begin{eqnarray*}
    d_H (\lbrack a,b\rbrack,\lbrack a_n,b_n \rbrack) &=& \max \left\lbrace
    \max_{x\in \lbrace a,b\rbrace} \min_{y\in \lbrace a_n,b_n\rbrace} \vert x - y \vert ,
    \max_{y\in \lbrace a_n,b_n\rbrace} \min_{x\in \lbrace a,b\rbrace} \vert x - y \vert \right\rbrace \\
    &=& \max \left\lbrace
    \vert a_n -a \vert,\vert b_n-b\vert \right\rbrace .
\end{eqnarray*}
(i)  We now construct the sequence $\Theta_n$ by letting $\Theta_n
= (a_n,b_n)$ where
\begin{equation*}
\begin{array}{cc}
a_n := \max \lbrace x\in p_n\, :\, x\leq a+ \delta_0\rbrace & b_n := \min \lbrace x\in p_n\, :\, x\geq b -
\delta_0\rbrace .
\end{array}
\end{equation*}
By the assumption made above, both sets are non - empty and by
Hausdorff - convergence we have $a_n\to a$, $b_n\to b$. Thus, by
(ii), $\overline{\Theta}_n$ converges to $\overline{\Theta}$.
(iii) By construction of the sequence in (i), always $y_n\leq
y_n^{\prime}$. But by (ii), $y_n,y_n^{\prime}\to y$ and
\begin{equation*}
    \vert y_n -  y_n^{\prime}\vert \leq \vert y_n -  y\vert + \vert y -  y_n^{\prime}\vert \leq 2 d_H
    (p_n,p).
\end{equation*}
\end{Proof}
As a first application, we derive that the counting function is lower semi-continuous.

\begin{Corollary}\label{cardinality} The function $\vert - \vert : \PP (U)\to\N$ is
lower semi-continuous if $\PP(U)$ is equipped with the Hausdorff topology.
\end{Corollary}

\begin{Proof} The intervals $\Theta_n\in\iota (p_n)$ approximating a given interval $\Theta\in\iota (p)$
constructed in Lemma \ref{convchar} can be chosen disjoint for
different $\Theta$. The left boundary points of the approximating
intervals are elements of $p_n$ and all different. Thus $\vert
p_n\vert\geq \vert p\vert$.
\end{Proof}

\subsection{Functions of Special Bounded Variation}\label{SBV}

Now we consider a special segmentation class, namely functions
which are of special bounded variation. They are defined as
follows: Recall that a right-continuous function of bounded
variation $f\in \BV (U)$ with left limits defines a signed measure
$\nu \in\MM (U)$ which is uniquely determined by it values $\nu
(\lbrack a,b)) := f(b) - f(a)$ on half-open intervals. Recall
further, that by {\em Lebesgue decomposition} (see for instance
\cite{Reed:Simon1:1980}, Theorem I.13, I.14, p. 22) there are
three uniquely determined measures $\nu_S,\nu_{R},\nu^{\perp}$
where $\nu = \nu_S+\nu_{R}+\nu^{\perp}$ and
\begin{enumerate}
\item $\nu_S =\sum c_x\,\delta_x$ is a sum of point measures,
\item $\nu_R << \lambda$ is absolutely continuous with respect to
Lebesgue measure with {\em Radon-Nikodym derivative} $f^{\prime} =
d\nu_R /d\lambda$, \item $\nu^{\perp}\perp\lambda$ is singular to
Lebesgue measure without point measures, i.e. $\nu^{\perp}
(\lbrace x \rbrace)=0$ for all $x\in U$.
\end{enumerate}

\begin{Definition} (1) A function $f\in \BV (U)$ is called to be of {\em
special bounded variation}, i.e. $f\in \SBV (U)$, if in the
decomposition above we have $\nu^{\perp} = 0$. It is called to be
of {\em p-special bounded variation}, i.e. $f\in\SBV_p(U)$, if in
addition, $f^{\prime}\in L^p (U)$ and -- in contrast to the usual
convention -- if the support $\supp \,(\nu_S)$ is a finite set.
(2) The {\em partition} associated to a function $f\in\SBV_2 (U)$
is given by $p(f):=\supp \,(\nu_S)\cup \lbrace 0,1\rbrace$.
\end{Definition}

It is not yet obvious that these functions form indeed a
segmentation class in the sense that they are equivalent to
piecewise continuous functions as assumed in Definition
\ref{segclass}. We will show this for $p=2$ and will characterize
them as piecewise Sobolev-functions.

\begin{Lemma} (1) Let $f\in\SBV_2(U)$. Then there are (up to an additive constant) uniquely determined functions
$t\in\TT (U)$, the space of step-functions introduced above, and
$F\in H^1 (\Omega)$, the Hilbert-Sobolev space of one time
generalized differentiable functions with square integrable
derivative on the open interval $\Omega := (0,1)$, such that
$f=t+F$. (2) In the equivalence class of $F\in H^1(\Omega)$, there
is a continuous representative as well denoted by $F\in C (U)$.
(3) The partition $p(f)$ associated to $f\in\SBV_2 (U)$ coincides
with $\disc (t)\cup \lbrace 0,1\rbrace$, where $\disc (t)$ denotes
the points of discontinuity of the corresponding step function.
\end{Lemma}

\begin{Proof} (1) Let $Df$ denote the distributional derivative (measure) of $f$.
By assumption $Df = f^{\prime}\cdot\lambda + \nu_S$. Then the
distribution function of the singular part $\nu_S$ is a right
continuous step function with left limits $t\in \TT (U)$. Hence
$D(f-t) = f^{\prime}\cdot\lambda$ with $\int_U \vert
f^{\prime}\vert^2 ds < \infty$. Hence $F:=f-t\in H^1 (\Omega)$.
(2) The second statement follows from {\em Sobolev's embedding
lemma} (see \cite{Reed:Simon2:1975}, Thm. IX.24, p. 52). (3) This
follows from the fact that the step function is the distribution
function of the point measure.
\end{Proof}

Thus, functions of special bounded variation are indeed
right-continuous with left limits and hence form a segmentation
class.

\section{The Reduction Principle}\label{RP}

The reduction principle consists of the fact that both, the
Mumford - Shah and the Blake - Zisserman functional -- except for
the degenerate case $\gamma = 0$ where this property is still true
in a restricted sense -- have the following property:

\vspace{0.1cm}

\noindent{\small\em Among all segmentations associated to a fixed
partition of the interval, the functional assumes a unique minimum.
The minimizer for a fixed partition can be computed by separately
minimizing independent functionals associated to the intervals of
the given partition.}

\vspace{0.1cm}

\noindent Hence -- by the first property -- the problem of
minimizing the whole functional can be reduced to the problem of
minimizing a {\em reduced functional} which is a function of
partitions rather than segmentations. Then, simple a priori bounds
on the number of jumps can be used to show that this partition
space is essentially compact. That, in particular, provides the
existence of global minimizers. Considerations like that are the
subject of this paper.

The second property is more important from the algorithmic point of
view. For discrete time-series and the {\em Potts model} this
property was used to establish an efficient algorithm to compute the
minimizers, see e.g. \cite{WiLi:2002} or the PhD-Thesis
\cite{Kempe:2004}. In the very recent PhD-Thesis
\cite{Friedrich:2004}, this property is used to construct efficient
algorithms in 2D when the admissible partitions are restricted to
certain subclasses (cf. the formulation of the reduction principle
in \cite{Friedrich:2004}, Definition 1.2.1).

\subsection{The Reduction Principle for Mumford - Shah}

We start with the statement of the reduction principle for Mumford
- Shah. First of all, we give the exact definition of the
functional already discussed above extending it to a functional on
$L^2 (U)$. Recall that functions of special bounded variation
provide a segmentation class as explained in Section \ref{SBV}.

\begin{Definition}\label{DefMS} Let $g\in L^2 (U)$, $\gamma,\mu \geq 0$. The {\em Mumford - Shah
functional} $\MS_{\gamma,\mu,g}: L^2 (U)\to\R$ to {\em signal} $g$
and parameters $\gamma,\mu$ is given by
\begin{equation}\label{MSdef}
    \MS_{\gamma,\mu,g}(f) := \left\lbrace\begin{array}{ll} \gamma j(f) +
    \frac{1}{\mu^2}\int_U \vert f^{\prime}\vert^2 ds + \Vert f -
    g\Vert^2 & f\in\SBV_2 (U) \\ \infty & \mathrm{else} \end{array}\right.
\end{equation}
where $j(f) := \vert p(f)\vert - 2$ is the {\em number of jumps}.
\end{Definition}

As proper segmentations of the signal $g$, we consider minimizers
of the functional, i.e.
\begin{equation*}
    f^{*}_{\gamma,\mu,\infty,g} := \argmin_{f\in L^2 (U)} \MS_{\gamma,\mu,g}(f).
\end{equation*}
The minimizer is not necessarily unique. The starting point for
the reduction principle is the fact that we may split up the
minimization procedure by the following observation: Let
\begin{equation}\label{observation}
    f^{*}_{\mu,\infty,g}(p) = \argmin_{\lbrace f\,:\, p(f)\subset p\rbrace}
\left\lbrack \frac{1}{\mu^2}\int_U \vert f^{\prime}\vert^2 ds +
\Vert f - g\Vert^2\right\rbrack .
\end{equation}
The key point, however, is that this minimizer for a fixed
partition exists, is unique and can be computed explicitly due to
a decoupling of the minimization procedure for the different
intervals in the complement of the partition. Then, the global
minimizer is given by
\begin{equation}\label{observationb}
    f^{*}_{\gamma,\mu,\infty,g} = \argmin_{f\in\lbrace f^{*}_{\mu,\infty,g}(p)\,:\,p\in\PP (U)\rbrace}\MS_{\gamma,\mu,g}(f) .
\end{equation}
That means, the minimization of the functional can be reduced to
the minimization on a much smaller subspace of the space of all
segmentations. This subspace is an image of the space of
partitions under the (in general not injective) map $p\mapsto
f^{*}_{\mu,\infty,g}(p)$.

First of all, we compute the unique minimizer in equation
(\ref{observation}).

\begin{Proposition}\label{reductionMS} Let $\gamma,\mu\geq 0$, $p\in\PP (U)$ be fixed. Then the unique minimizer
\begin{equation*}
    f_{\mu,\infty,g}^{*} (p) := \argmin_{\lbrace f\,:\, p(f)\subset p\rbrace}
\left\lbrack \frac{1}{\mu^2}\int_U \vert f^{\prime}\vert^2 ds +
\Vert f - g\Vert^2\right\rbrack
\end{equation*}
can be constructed as follows: Let $\Theta\in\iota (p)$. Denote by
$\Delta_{\Theta}$ the Laplacian on $\Theta$ with {\em Neumann
boundary conditions}. Then the function
\begin{equation*}
    f_{\Theta} := - \mu^2 R (\Delta_{\Theta},\mu^2) g
\end{equation*}
where $R(\Delta_{\Theta},\mu^2)$ denotes the {\em resolvent} of the Laplacian, is continuous on $\overline{\Theta}$.
Denote by $\Theta^+$, the closed interval $\overline{\Theta}$ with right boundary point removed. Then
\begin{equation}\label{miniMS}
    f_{\mu,\infty,g}^{*} (p) = \sum_{\Theta\in \iota (p)} f_{\Theta} \chi_{\Theta^+} .
\end{equation}
\end{Proposition}

\begin{Remark} (i) Note that $p(f_{\mu,\infty,g}^{*}
(p^{\prime}))\subset p^{\prime}$ but that $p(f_{\mu,\infty,g}^{*}
(p^{\prime}))$ can be strictly smaller. That explains why the map
$p\mapsto f_{\mu,\infty,g}^{*} (p)$ is not injective and was the
reason to put $\lbrace f\,:\, p(f)\subset p\rbrace$ in equation
(\ref{observation}). (ii) Having in mind the decomposition
$L^2(U)=\oplus_{\Theta\in\iota(p)} L^2 (\Theta^+)$ we could write
shortly
\begin{equation*}
f_{\mu,\infty,g}^{*} (p) = -\mu^2 \bigoplus_{\Theta\in \iota (p)}
R(\Delta_{\Theta},\mu^2)g
\end{equation*}
where we identify $R(\Delta_{\Theta},\mu^2)g\in H^1 (\Theta)$ with
its image in $C(\Theta^+)$ according to the Sobolev embedding
theorem. The second property of the reduction principle is
reflected by the fact that the operator assigning to the signal
$g$ the minimizer $f_{\mu,\infty,g}^{*} (p)$ is {\em decomposable}
(see for instance \cite{Reed:Simon4:1978}, p. 284 ff.) with
respect to the orthogonal sum decomposition
$L^2(U)=\oplus_{\Theta\in\iota(p)} L^2 (\Theta^+)$. (3) Note that
the case $\mu = 0$ is included in the preceding considerations in
the following sense: If $\mu = 0$, the penalization for
$\SBV$-functions with non - vanishing derivative tends to
infinity, hence constant functions $f_{\Theta}$ are favored in
this case. According to this, the operator $-\mu^2
\,R(\Delta_{\Theta},\mu^2)$ tends to the projection onto the
kernel of $\Delta_{\Theta}$ which consists of all constant
functions on $\Theta$.
\end{Remark}

\begin{Proof} Minimizing the Mumford - Shah functional for a fixed partition means minimization of the expression
\begin{equation*}
\sum_{\Theta\in\iota (p)} \int_{\Theta} dx \left(\mu^{-2} \vert
f^{\prime}_{\Theta}\vert^2 + \vert f_{\Theta} - g\vert^2 \right)
\end{equation*}
where we first consider $f_{\Theta,\Theta\in\iota (p)}$ to be a
tuple of functions $f_{\Theta}\in C^{\infty}(\overline{\Theta})$.
Let $h_{\Theta,\Theta\in\iota (p)}$ be another such tuple and
$\tau_{\Theta , \Theta\in\iota (p)} \in \R^{\vert\iota (p)\vert}$.
The minimization condition reads
\begin{eqnarray*}
0 &=& \frac{\partial}{\partial \tau_{\Theta}}\left. \int_{\Theta} dx \left(\mu^{-2} \vert f^{\prime}_{\Theta} +
\tau_{\Theta} h^{\prime}_{\Theta}\vert^2 + \vert f_{\Theta} + \tau_{\Theta} h_{\Theta} -
g\vert^2 \right)\right\vert_{\tau_{\Theta}=0} \\
&=& 2 \int_{\Theta} dx \left(\mu^{-2}  f^{\prime}_{\Theta}
h^{\prime}_{\Theta} +  f_{\Theta} - g \right) \\
&=& 2 \left( -\mu^{-2}\langle f^{\prime}_{\Theta},
h_{\Theta}^{\prime}\rangle + \langle f_{\Theta}, h_{\Theta}\rangle
- \langle g ,h_{\Theta}\rangle \right) .
\end{eqnarray*}
Hence, adopting the notation from \cite{Taylor:1996}, Section 5.7,
p. 345 ff, we may consider the extended map $-\mu^2 \mathcal{L}_N
+ 1 $ from $H^1 (\Theta)$ to the dual space $H^1 (\Theta)^*$
defined by the relation
\begin{equation*}
    \langle (-\mu^2 \mathcal{L}_N
+ 1) f_{\Theta}, h_{\Theta}\rangle := \mu^{-2}\langle
    f_{\Theta}^{\prime}, h_{\Theta}^{\prime}\rangle + \langle
    f_{\Theta}, h_{\Theta}\rangle .
\end{equation*}
By \cite{Taylor:1996}, Proposition 7.2, p. 346, for $\mu > 0$,
$g\in L^2 (U)$, the equation is solved by a unique $f_{\Theta}\in
H^2 (\Theta)$ satisfying
\begin{equation*}
\begin{array}{lr}
-\mu^2 f^{\prime\prime}_{\Theta} + f_{\Theta} = g & \mathrm{on}\,\,\Theta \\
f^{\prime}_{\Theta}\vert_{\partial\Theta} = 0 & \\
\end{array}
\end{equation*}
where the second equation is understood in terms of the {\em trace
map} (see \cite{Taylor:1996}, Proposition 1.6, p. 273. Thus, the
Euler equation can be equivalently described by
$(\Delta_{\Theta}-\mu^2)f_{\Theta} = -\mu^2 g$ where
$\Delta_{\Theta}$ is the {\em Neumann laplacian} on $\Theta$.
Hence
\begin{equation*}
    f_{\Theta} = -\mu^2 (\Delta_{\Theta} - \mu^2 )^{-1} g = - \mu^2 R(\Delta_{\Theta} , \mu^2)
    g
\end{equation*}
where $R$ denotes the {\em resolvent}. By {\em Sobolev's embedding theorem},
\begin{equation*}
f_{\Theta}\in\dom (\Delta_{\Theta}) \subseteq  H^1 (\Theta)
\subset C (\overline{\Theta}),
\end{equation*}
and $f^{*}$ provides indeed a right-continuous solution. In the
case $\mu = 0$, the minimum can only be assumed by locally
constant functions, i.e. $f_{\Theta}^{\prime}=0$. Hence, in that
case the only non-trivial variations of $f_{\Theta}$ are given by
constant functions $h_{\Theta}$ as well. Thus, the solution to the
variational problem is given by the constant function assuming the
mean value of $g$ on $\Theta$, i.e.
\begin{equation*}
    \langle f_{\Theta}, h_{\Theta}\rangle = f_{\Theta}h_{\Theta}\langle
    1,1\rangle = \langle
    g,h_{\Theta}\rangle = h_{\Theta}\langle
    g,1\rangle
\end{equation*}
and therefore $f_{\Theta}= \langle g,1\rangle/\langle 1,1\rangle$.
The minimizer depends thus even continuously from $\mu$, since the
resolvent function considered above tends to the projection onto
the kernel of the Dirichlet Laplacian as $\mu$ tends to zero.
Since the kernel consists of constant functions, this coincides
with the mean value.
\end{Proof}

Hence, the precise formulation of the reduction principle in the
case of the Mumford - Shah functional is given by

\begin{Corollary}[Reduction for Mumford-Shah]\label{RPMS} Let $\gamma,\mu \geq 0$. The minimization of
the Mumford - Shah functional is equivalent to the minimization of
the {\em reduced Mumford - Shah functional} $ \ms_{\gamma,\mu,g} :
\PP (U) \to \R$ given by
\begin{equation*}
    \ms_{\gamma,\mu,g} (p):=\MS_{\gamma,\mu,g}(f_{\mu,\infty,g}^{*}
    (p))= \gamma \,j(p)  - \langle u , f_{\mu,\infty,g}^{*}
    (p)\rangle + \Vert g \Vert^2 ,
\end{equation*}
more precisely, $p$ is a minimizer of $\ms_{\gamma,\mu,g}$ if and
only if $f_{\mu,\infty,g}^{*}(p)$ is a minimizer of
$\MS_{\gamma,\mu,g}$.
\end{Corollary}

\subsection{The Reduction Principle for Blake - Zisserman}

The reduction principle for the Blake - Zisserman functional is
similar to the one for Mumford - Shah. We just have to identify
the respective quantities in the discrete setting. Again, we start
by giving our definition of the functional. Recall the notions of
{\em equidistant setup} and {\em conditional expectation} from
Definition \ref{SetUp} and the classes of considered
step-functions from Definition \ref{Treppen}.

\begin{Definition} \label{BZD} Let $n\in\N$, $g\in L^2 (U)$ and $\gamma,\mu \geq 0$ be
fixed. The {\em Blake-Zisserman functional} $\BZ_{\gamma,\mu,g}^n
: L^2 (U) \to \R$ is given by
\begin{equation}\label{BZdef}
\BZ_{\gamma,\mu,g}^n (f) := \left\lbrace\begin{array}{ll}
\Phi_{\gamma,\mu}^n (f) + \frac{1}{n}\sum_{\kappa =0}^{n-1} (f -
g_n )^2(\kn) & f\in \TT_n(U) \\ \infty &\mathrm{else}
\end{array} \right.
\end{equation}
where
\begin{equation*}
    \Phi_{\gamma,\mu}^n (f) := \sum_{\kappa =0}^{n-2} \min\lbrace \koeff(f(\kpn) -
    f(\kn))^2,\gamma\rbrace
\end{equation*}
and $g_n = E(g\,\vert\,\sigma_n)$ denotes conditional expectation.
\end{Definition}

Again, we seek for minimizers. The starting point for the
reduction principle in this case will as well consist of an
observation concerning the minimization of the functional. Let
\begin{equation*}
    f^{*}_{\gamma,\mu,n,g} := \argmin \BZ_{\gamma,\mu,g}^n (f)
\end{equation*}
be a (not necessarily unique) minimizer. Then we have

\begin{Lemma}\label{edgemod} The minimization of $\BZ_{\gamma,\mu,g}^n$ is
equivalent to the minimization of the functional
$H_{\gamma,\mu,g}^n : \Z_2^{n-1}\times \TT_n (U)\to \R$ given by
\begin{equation*}
    H_{\gamma,\mu}^n (e,f) := \Psi_{\gamma,\mu}^n (e,f)+ \frac{1}{n}\sum_{\kappa =0}^{n-1} \koeff
    (f-g_n)^2(\kn)
\end{equation*}
where $\Z_2 :=\lbrace 0,1\rbrace$, $e= (e_1,...,e_{n-1})$ and
\begin{equation*}
    \Psi_{\gamma,\mu}^n (e,f):=\sum_{\kappa =0}^{n-2} \koeff (f(\kpn)-f (\kn))^2
    (1-e_{\kappa +1}) + \gamma\, e_{\kappa +1}.
\end{equation*}
\end{Lemma}

\begin{Proof} See \cite{Winkler:2002}, p. 36 f.
\end{Proof}

The points $\kn$ such that $e_{\kappa} = 1$ correspond to {\em
edges} in the segmentation of $u_n$ (see again the presentation in
\cite{Winkler:2002}, p. 36 f). They correspond exactly to those
points $\kn$, where the minimizer fulfills $\vert
f^{*}_{\gamma,\mu,n,g}(\kpn)-f^{*}_{\gamma,\mu,n,g} (\kn)\vert
>\mu\,\sqrt{\gamma /n}$, i.e. they represent the {\em discrete discontinuities} of $f^{*}_{\gamma,\mu,n,g}$ with respect to the threshold
$T\equiv\mu\,\sqrt{\gamma/n}$ as defined in Definition
\ref{discdef}. Furthermore, there is a bijection onto associated
partitions
\begin{equation*}
    p(e) := \lbrace \kn \,:\,e_{\kappa} = 1\rbrace \cup\lbrace 0,1\rbrace\in\PP_n (U)
\end{equation*}
and these form the discrete analogue of the partitions in the
continuous case. Using the lemma above, the starting point for the
reduction principle is provided as in the continuous case by the
following observation: Let
\begin{equation*}
    Q^n_{\mu}(e,f) := \sum_{\kappa =0}^{n-2} \koeff (f^{\kappa +1} - f^{\kappa})^2 (1-e_{\kappa})
    + \frac{1}{n}
    \sum_{\kappa =0}^{n-1} (f^{\kappa} - g_n (\kn))^2
\end{equation*}
and
\begin{equation}\label{observation2}
    f^{*}_{\mu,n,g} (p) = \argmin_{\lbrace f :
    p_{\sqrt{\gamma\mu /n}}(f) \subset p\rbrace}
    Q^n_{\mu}(e,f).
\end{equation}
Again, the minimizer for a fixed partition exists, is unique and
can be computed explicitly due to a decoupling of the minimization
procedure for the different intervals in the complement of the
partition. Then, the global minimizer is given by
\begin{equation}\label{observation2b}
    f^{*}_{\gamma,\mu,n,g} = \argmin_{f\in\lbrace f^{*}_{\mu,n,g}(p)\,:\,p\in\PP_n (U)\rbrace}\BZ^n_{\gamma,\mu,g}(f) .
\end{equation}

First of all, we recall the decomposition of the conditional
expectation map that was already introduced in Section \ref{1.2}.

\begin{Definition}\label{pichi} The map $\pi_n : L^2 (U)\to \R^n$
is given by the row
\begin{equation*}
    \pi_n g := (g_n(0), g_n (1/n),..., g_n (n-1/n))
\end{equation*}
where $g_n := E(g\,\vert\,\sigma_n)$ denotes conditional
expectation. The map $\tau_n : \R^n \to L^2 (U)$ is given by
\begin{equation*}
    \tau_n{\underline{f}}(s) := f^{n-1}\chi_{\RI}(s) + \sum_{\kappa =0}^{n-2} f^{\kappa} \,\chi_{\lbrack\kn,\kpn )}(s).
\end{equation*}
\end{Definition}

\begin{Remark} Clearly, $g_n = \tau_n\circ\pi_n g$.
\end{Remark}

As in the continuous case, we start by constructing the unique
minimizer associated to a fixed partition.

\begin{Proposition}\label{reductionBZ} Let $\gamma,\mu \geq 0$ and $p\in\PP_n (U)$ be fixed. Then the minimizer
\begin{equation*}
    f^{*}_{\mu,n,g} (p) = \argmin_{\lbrace f :
    p_{\sqrt{\gamma\mu /n}}(f) \subset p\rbrace} Q^n_{\mu}(e,f)
\end{equation*}
is unique and can be constructed as follows: Let $p = \lbrace 0 =:
\kappa_0/n < \kappa_1/n < \kappa_2/n < ... < \kappa_l/n < 1 =:
k_{l+1}/n \rbrace\in \PP_n(U)$. Consider the block matrix
\begin{equation*}
    A(p) := \left( \begin{array}{cccc} B(1) & & 0 \\ & \ddots &  \\0  & & B(l+1)\\\end{array}\right)
\end{equation*}
where the blocks $B(r)$ are given by the $\kappa_r -
\kappa_{r-1}\times \kappa_r - \kappa_{r-1}$- band matrices
\begin{equation*}
    B(r) := \left( \begin{array}{cccccc} -1 & 1 & & &  0\\
    1 & -2 & 1 & &   \\ & \ddots & \ddots & \ddots &   \\
     & & 1 & -2 & 1\\ 0 & & & 1 & -1\\\end{array}\right) .
\end{equation*}
Then
\begin{equation*}
    f_{\mu,n,g}^{*} (p)  := - \mu^2 \tau_n\circ R (n^2A(p),\mu^2)
    \circ\pi_n g^t
\end{equation*}
where $R (n^2A(p),\mu^2)$ denotes the resolvent and $\pi_n g^t$
the transpose vector.
\end{Proposition}

\begin{Remark} As in the continuous case, the block structure of
$A(p)$ corresponds to a direct sum decomposition of underlying
space and operator.
\end{Remark}

\begin{Proof} Analogous to the proof of Proposition \ref{reductionMS} minimization of $\BZ$
for fixed partition $e$, writing shortly $f^{\kappa} = f (\kn)$,
equivalent to
\begin{eqnarray*}
    0 &=& \frac{\partial}{\partial f^{\kappa}} \sum_{\kappa =0}^{n-2} \koeff (f^{\kappa +1} - f^{\kappa})^2 (1-e_{\kappa})
    + \gamma e_{\kappa} + \frac{1}{n}
    \sum_{\kappa =0}^{n-1} (f^{\kappa} - g_n (\kn))^2 \\
    &=& \koeff 2((f^{\kappa +1} - f^{\kappa})(1-e_{\kappa}) + (f^{\kappa} - f^{\kappa -1})(1-e_{\kappa -1})) + \frac{2}{n}
     (f^{\kappa} - g_n (\kn))
\end{eqnarray*}
for all $\kappa =0,...,n-1$. This system of linear equations can
be written as
\begin{equation*}
    (n^2 A(p) - \mu^2) \underline{f} = -\mu^2 \pi_ng^t,
\end{equation*}
hence
\begin{equation*}
     \underline{f}^*_{\mu,n,g} = \tau_n (\underline{f})=-\mu^2 \tau_n\,(n^2 A(p) - \mu^2)^{-1}\pi_n g^t =
     -\mu^2 \tau_n\,R(n^2 A(p), \mu^2)\,\pi_n g^t.
\end{equation*}
\end{Proof}

Hence, we obtain in analogy to Corollary \ref{RPMS}

\begin{Corollary}[Reduction for Blake-Zisserman]\label{RPBZ} \hspace{0.2cm} Let $\gamma,\mu \geq 0$. The minimization of the Blake - Zisserman functional is equivalent
to the minimization of the {\em reduced Blake - Zisserman
functional} $ \bz_{\gamma,\mu,g}^n : \PP (U) \to \R$ given by
\begin{equation*}
    \bz_{\gamma,\mu,g}^n (p):=\left\lbrace\begin{array}{ll}\BZ_{\gamma,\mu,u}^n(f_{\mu,n,g}^{*}
    (p)) & p\in\PP_n (U) \\ \infty & \mathrm{else}\end{array}\right..
\end{equation*}
For $p\in\PP_n (U)$, the reduced functional is given by
\begin{equation*}
    \bz_{\gamma,\mu,g}^n (p) = \gamma\, j(p) - \langle
    g_n , f_{\mu,n,g}^{*}
    (p)\rangle + \Vert g_n\Vert^2.
\end{equation*}
\end{Corollary}

\section{Continuity of the Reduced Functionals}

So far, we have seen that the minimization on the function space
can be reduced to a minimization on partition space. However, it
is not yet clear whether the reduced functionals $\ms$ and $\bz$
depend continuously on the parameters and/or the partitions. In
order to prove this, we first have to investigate the continuity
properties of the minimizers for a fixed partition -- from now on
shortly denoted {\em partition solvers} -- introduced in the
preceding section.

\subsection{The Partition Solvers}

We consider continuity properties of the partition solvers. It
turns out that they depend continuously on the parameters. Since
the partition solvers for the Blake - Zisserman functional can not
be applied to all partitions, we first have to define a proper
domain.

\begin{Definition}\label{partitiondomain}  Recall the definition
of
\begin{equation*}
    \T:= \lbrace 1/n\,:\,n\in\N\rbrace\cup\lbrace 0 \rbrace \subset\R
\end{equation*}
from Theorem \ref{Hauptsatz}. Adopting the convention that
$\PP_{1/0}(U) = \PP (U)$ let
\begin{equation*}
    \E(\T):= \lbrace (t,p)\in \T\times \PP(U)\,:\, p\in
    \PP_{1/t}(U)\rbrace\subset \R\times\PP(U)
\end{equation*}
equipped with its relative topology, i.e. all accumulation points
are of the form $(0,p)$, $p\in\PP(U)$.
\end{Definition}

With these conventions, the result about the parameter dependence
of the partition solvers reads as follows.

\begin{Theorem}\label{maineins} For fixed $g\in L^2 (U)$, the {\em partition solver map} $f^{*} :
\R^+_0\times  \E (\T)\rightarrow L^2 (U)$, given by
\begin{equation*}
     f^{*} (\mu,t,p) := f^{*}_{\mu , 1/t,g}(p)
\end{equation*}
is continuous.
\end{Theorem}

\begin{Proof} The statement follows from the series of Lemmas proved below.
By Lemma \ref{psolverzwei}, \ref{psolverdrei}, $f^{*}$ is
continuous in $\E (\T)$ for fixed $\mu$. By Lemma
\ref{psolvereins}, $f^{*}$ is equicontinuous in $\mu$ for all
converging sequences of partitions. That implies joint continuity.
\end{Proof}

Before we come to the lemmas implyint Theorem \ref{maineins}, we
apply the theorem to prove lower semi - continuity of the reduced
Mumford - Shah functional.

\begin{Corollary}\label{semicont} The reduced functional $\ms_{\gamma,\mu,g}$ is lower semi - continuous.
\end{Corollary}

\begin{Proof} By Corollary \ref{cardinality}, the function
$j:\PP(U)\to \N$ is lower semi - continuous. By Theorem
\ref{maineins}, $f^{*}$ depends for fixed $\gamma,\mu,t=0$
continuously on $p\in\PP(U)$. By the explicit form of
$\ms_{\gamma,\mu,g}$ given in Corollary \ref{RPMS}, that implies
the statement.
\end{Proof}

The first lemma states that a family of partition solvers
associated to a convergent sequence of partitions depends
equicontinuously on the parameter $\mu\geq 0$.

\begin{Lemma}\label{psolvereins} (i) Let $p\in\PP (U)$ be fixed. Then the map $G_{\infty}:\R^+_0\to  L^2 (U)$ given by
$G_{\infty}(\mu):= f^{*}_{\mu,\infty,g} (p)$ is Lipschitz
continuous. (ii) Let $p\in\PP_n(U)$ be fixed. Then the map $G_n
(\mu) := f^{*}_{\mu,n,g} (p): \R^+_0 \to  L^2 (U)$ is Lipschitz
continuous. (iii) Let $p_{n,n\in\N}$ be a family of partitions
such that $p_n\in\PP_n(U)$ converges to $p\in\PP(U)$ with respect
to Hausdorff metric. Then the associated family
$G_{n,n=1,2,...,\infty}$ is uniformly Lipschitz.
\end{Lemma}

\begin{Proof} (i) Recall the notations from Proposition \ref{reductionMS}. By definition, all partitions are finite. Thus
it is enough to show Lipschitz continuity for one single interval,
i.e. for a subinterval $\Theta := (x,y)\in\iota (p)$, the map
\begin{equation*}
    f_{\Theta} (\mu) := -\mu^2 R (\Delta_{\Theta},\mu^2)g
\end{equation*}
is Lipschitz continuous. Since the {\em Neumann laplacian} is
self-adjoint and has a discrete spectrum with semisimple
non-positive eigenvalues and its kernel consists exactly of
constant functions, we have the spectral decomposition
\begin{equation*}
    -\mu^2 R (\Delta_{\Theta},\mu^2)g = E_0 g + \sum_{s\geq 1}\frac{\mu^2}{\mu^2 -
    \lambda_s(\Theta)}E_s g
\end{equation*}
where $E_0 g$ denotes the {\em mean value} of $g$ in $\Omega$,
i.e. the orthogonal projection of $g$ onto constant functions.
Hence, by $0
 > \lambda_1 > \lambda_2 ...$, we obtain
\begin{eqnarray*}
    \Vert f_{\Theta} (\mu) - f_{\Theta} (\mu^{\prime}) \Vert^2
    &=& \left\Vert\sum_{s\geq 1}\left\lbrack\frac{\mu^2}{\mu^2 -
    \lambda_s({\Theta})} - \frac{\mu^{\prime 2}}{\mu^{\prime 2} -
    \lambda_s({\Theta})}\right\rbrack E_s g\right\Vert^2 \\
    &=& \sum_{s\geq 1}\frac{\lambda_s^2({\Theta})(\mu^{\prime 2}-\mu^2)^2}{(\mu^2\mu^{\prime 2}
    - \lambda_s({\Theta})(\mu^2 + \mu^{\prime 2}) + \lambda_s^2(\Theta))^2}  \left\Vert E_s g\right\Vert^2 \\
    &\leq& \sum_{s\geq 1}\frac{(\mu^{\prime 2}-\mu^2)^2}{\lambda_s^2(\Theta)}  \left\Vert E_s g\right\Vert^2
    \leq \frac{(\mu^{\prime 2}-\mu^2)^2}{\lambda_1^2(\Theta)}  \Vert g\Vert^2 .
\end{eqnarray*}
Hence $f_{\Theta}$ is Lipschitz continuous with Lipschitz constant
$C_{\Theta}:=\Vert g\Vert /\vert\lambda_1 (\Theta)\vert$ depending
on the norm of $g$ and on the {\em spectral gap} $\vert\lambda_1
(\Theta)\vert$ of the Neumann Laplacian on $\Theta$. That implies
\begin{equation*}
    \Vert G_{\infty}(\mu) - G_{\infty}(\mu^{\prime})\Vert \leq
    C \vert \mu^2 - \mu^{\prime 2}\vert
\end{equation*}
where $C:= \max_{\Theta\in\iota (p)} C_{\Theta}<\infty $. (ii)
Recall the notation from Proposition \ref{reductionBZ}. By
definition, all partitions are finite. Thus it is enough to show
Lipschitz continuity for one single interval, i.e. for a
subinterval $\Theta_r := \lbrack
\kappa_{r-1}/n,\kappa_{r}/n)\subset U$ and the corresponding map
\begin{equation*}
    f_{\Theta_r} (\mu) :=  -\mu^2 \tau_n \,R (n^2
    B(r),\mu^2)\,\pi_n g^t.
\end{equation*}
$B(r)$ is a symmetric matrix and its kernel consists exactly of
constant vectors. Thus, we have the spectral decomposition
\begin{equation*}
    -\mu^2 R (n^2
    B(r),\mu^2)\pi_n g^t = P_0 \pi_n g^t + \sum_{s = 1}^{\kappa_{r-1} - \kappa_{r -1}-1}\frac{\mu^2}{\mu^2 -
    n^2\lambda_s(B(r))}P^n_s \pi_n g^t
\end{equation*}
where $P_0 u$ denotes the orthogonal projection of
$\underline{u}_n$ onto constant vectors. Hence we obtain using the
fact that conditional expectation is a projection (which implies
$\Vert g_n\Vert \leq \Vert g \Vert$)
\begin{eqnarray*}
    & & \Vert f_{\Theta_r} (\mu) - f_{\Theta_r} (\mu^{\prime}) \Vert^2 \\
    &=&
    \sum_{s = 1}^{\kappa_{r+1} - k\kappa_r -1}\left\lbrack\frac{\mu^2}{\mu^2 -
    n^2\lambda_s(B(r))} - \frac{\mu^{\prime 2}}{\mu^{\prime 2} -
    n^2\lambda_s(B(r))}\right\rbrack^2 \Vert \tau_n P_s^n
    \pi_n g^t\Vert^2\\
    &\leq& \frac{(\mu^{\prime 2}-\mu^2)^2}{n^4\lambda_1^2(B(r)))}  \Vert g\Vert^2
\end{eqnarray*}
where $\vert\lambda_1 (B(r))\vert$ is the {\em spectral gap} of
$B(r)$. By the same argument as above, $G_n$ is Lipschitz
continuous. (iii) By the preceding parts of the proof and by Lemma
\ref{convchar}, it is enough to show the following:

\noindent{\small\em Let $\overline{\Theta}_n := \lbrack
\kappa_r^{(n)} /n,\kappa_{r+1}^{(n)}/n \rbrack$ be a sequence of
intervals tending to $\overline{\Theta} := \lbrack \kappa^-,
\kappa^+\rbrack$ in Hausdorff metric. Then the spectral gap
$n^2\,\vert\lambda_1 (B^{(n)}(r))\vert$ of the associated block
matrices is uniformly bounded below.}

Note that the case $\kappa^+=\kappa^-$, i.e.
$\overline{\Theta}=\lbrace\kappa^+\rbrace$, is included in the
following considerations. First of all, by Lemma \ref{convchar},
Hausdorff convergence $\overline{\Theta}_n\to \overline{\Theta}$
is equivalent to $\lim \kappa^{(n)}_r/n = \kappa^-$, $\lim
\kappa^{(n)}_{r+1}/n = \kappa^+$. That implies
\begin{equation*}
    \frac{\kappa^{(n)}_{r+1} - \kappa^{(n)}_r}{n} = \kappa^+ -
    \kappa^- + o(n^{-1}).
\end{equation*}
On the other hand, the eigenvalues of $B^{(n)}(r)$ are well known
to be
\begin{equation*}
    \lambda_s := 2\left(\cos \left({\scriptstyle \frac{\pi(s-1)}{\kappa_r - \kappa_{r-1}}}\right)
    -1\right)
\end{equation*}
where $s=1,...,\kappa_{r}- \kappa_{r-1}$ (see e.g.
\cite{Kuensch:94}, Theorem 1.3). Thus, the largest non-zero
eigenvalue is given by
\begin{equation*}
    \lambda_1 (B^{(n)}(r)) = 2\left(\cos \left({\scriptstyle \frac{\pi}{\kappa_r - \kappa_{r-1}}}\right)
    -1\right)
\end{equation*}
and hence
\begin{eqnarray*}
    n^2 \lambda_1 (B^{(n)}(r))
    &=&  2n^2\left(\cos ({\scriptstyle\frac{\pi}{\kappa_r - \kappa_{r-1}}})
    -1\right) =  2n^2\left(\cos ({\scriptstyle\frac{\pi}{n(\kappa^+ - \kappa_{-} +
    o(n^{-1}))}})
    -1\right) \\ &=&  - \frac{\pi^2}{(\kappa^+ - \kappa_{-})^2} + o(n^{-1}).
\end{eqnarray*}
That implies finally that for all $\epsilon > 0$ there is some
$n_0$ with
\begin{equation*}
    C_n \leq \frac{\Vert g \Vert}{\pi^2} L^2 + \epsilon
\end{equation*}
for all $n>n_0$ where $L$ is the maximal length of an interval in
$\iota (p)$. Hence, the family $G_n$ is uniformly Lipschitz.
\end{Proof}

The next result shows that the partition solver associated to the
Mumford - Shah functional depends continuously on the given
partition.

\begin{Lemma}\label{psolverzwei} The function
$f_{\mu,\infty,g}^{*}:\PP(U)\to \SBV_2 (U)\subset L^2(U)$ is
continuous if $\PP(U)$ is equipped with the Hausdorff topology.
\end{Lemma}

\begin{Proof} First of all, recall that there is a strongly
continuous representation of the {\em affine group}
$\mathrm{Aff}(\R) := \lbrace L_{a,b}(x):= ax + b\,:\, a>0,
b\in\R\rbrace$ on $L^2 (\R)$ given by $\rho (L_{a,b})f := f
(ax+b)$. In particular, given a sequence $(a_s,b_s)$ tending to
$(1,0)$ as $s$ tends to infinity, we obtain
\begin{equation*}
    \lim_{s\to\infty} \Vert \rho (L_{a_s,b_s})f - f\Vert = 0
\end{equation*}
for all $f\in L^2 (\R)$. By Lemma \ref{convchar}, it is again
sufficient to prove that for a sequence $\Theta_{n}:=
(\kappa^{(n)}_- , \kappa^{(n)}_+)$ of intervals such that
$\overline{\Theta}_n\to \overline{\Theta}:=\lbrack \kappa_- ,
\kappa_+\rbrack$ in Hausdorff distance, we have
\begin{equation*}
    \lim_{n\to\infty} \Vert \mu^2 R(\Delta_{\Theta},\mu^2)g -
    \mu^2 R(\Delta_{\Theta_n},\mu^2)g\Vert = 0.
\end{equation*}
(i) Let $\kappa_+\neq \kappa_-$ and therefore $L_n(x):= w_n x +
v_n :\R\to\R$ the unique linear map such that
$L_n(\kappa^{(n)}_-)=\kappa_-$ and $L_n(\kappa^{(n)}_+)=\kappa_+$.
Now $L_n$ maps the domains of the corresponding Neumann
laplacians, i.e. we have $g\circ L_n\in\DD (\Delta_{\Theta_n})
\Leftrightarrow g\in\DD (\Delta_{\Theta})$, and furthermore
\begin{equation*}
    \Delta_{\Theta_n} (g\circ L_n) = w_n^2\,
    \,(\Delta_{\Theta}g)\circ L_n .
\end{equation*}
In the spirit of the remark made above about the representation of
the affine group, we write $\rho_n (f):= f\circ L_n$. In order to
avoid difficulties with the domains, we will write $g_{\Theta}$
for the restriction $g\chi_{\Theta}$ of a function $g$ to the
interval $\Theta$ having in mind that due to the support of the
resolvent kernel, we always have that
\begin{equation*}
    R(\Delta_{\Theta},\mu^2)g = R(\Delta_{\Theta},\mu^2)g_{\Theta}
\end{equation*}
is supported on $\Theta$. Now
\begin{equation*}
    R(\Delta_{\Theta_n},\mu^2) (g_{\Theta}\circ L_n) = w_n^{-2}
    \,(R(\Delta_{\Theta},w_n^{-2}\mu^2)g_{\Theta})\circ L_n
\end{equation*}
and by $\Vert g_{\Theta}\circ L_n\Vert = w_n^{-1/2} \Vert
g_{\Theta} \Vert$ and $\Vert\mu^2 R(\Delta_{\Theta},\mu^2) \Vert =
1$ for all $\mu\geq 0$, $\Theta\in\iota (p)$ we therefore obtain
letting $g(n)_{\Theta} := g_{\Theta_n}\circ L_n^{-1}=
\rho_n^{-1}(g_{\Theta_n})$
\begin{eqnarray*}
    & &\Vert \mu^2 R(\Delta_{\Theta},\mu^2)g_{\Theta} -
    \mu^2 R(\Delta_{\Theta_n},\mu^2)g_{\Theta_n}\Vert \\
    &=& \Vert \mu^2 R(\Delta_{\Theta},\mu^2)g_{\Theta} -
    \left(\mu^2 w_n^{-2} R(\Delta_{\Theta},\mu^2 w_n^{-2})g(n)_{\Theta}\right)\circ L_n\Vert \\
    &\leq& \Vert \mu^2 R(\Delta_{\Theta},\mu^2)g_{\Theta} - (\mu^2 R(\Delta_{\Theta},\mu^2)g_{\Theta})\circ
    L_n\Vert \\ & & + \Vert (\mu^2 R(\Delta_{\Theta},\mu^2)(g_{\Theta} - g(n)_{\Theta}))\circ
    L_n\Vert \\& & + \left\Vert\left((\mu^2 R(\Delta_{\Theta},\mu^2) -
    \mu^2 w_n^{-2} R(\Delta_{\Theta},\mu^2 w_n^{-2}))g(n)_{\Theta}\right)\circ
    L_n\right\Vert\\
    &\leq& \Vert \mu^2 R(\Delta_{\Theta},\mu^2)g_{\Theta} - \rho_n(\mu^2 R(\Delta_{\Theta},\mu^2)g_{\Theta})\Vert \\
    & & + \Vert \rho_n\left(\mu^2 R(\Delta_{\Theta},\mu^2)(g_{\Theta} - g(n)_{\Theta})\right)\Vert \\
    & & + \left\Vert\rho_n\left((\mu^2 R(\Delta_{\Theta},\mu^2) -
    \mu^2 w_n^{-2} R(\Delta_{\Theta},\mu^2 w_n^{-2}))g(n)_{\Theta}\right)\right\Vert\\
    &\leq& \Vert \mu^2 R(\Delta_{\Theta},\mu^2)g_{\Theta} - \rho_n(\mu^2 R(\Delta_{\Theta},\mu^2)g_{\Theta})\Vert \\
    & & + \Vert \rho_n\Vert \,\Vert g_{\Theta} - g(n)_{\Theta}\Vert + 2\mu^2 \Vert\rho_n\Vert \,
    \Vert g(n)_{\Theta}\Vert\,\vert 1 -  w_n^{-2}\vert\\
\end{eqnarray*}
using the {\em resolvent identity} for the final step. Note that
\begin{eqnarray*}
\Vert g_{\Theta} - g(n)_{\Theta}\Vert &=& \Vert g_{\Theta} -
\rho_n^{-1}(g_{\Theta}) + \rho_n^{-1}(g_{\Theta})- \rho_n^{-1}
(g_{\Theta_n}) \Vert \\
&\leq & \Vert g_{\Theta} - \rho_n^{-1}(g_{\Theta}) \Vert + \Vert
\rho_n^{-1}(g_{\Theta}- g_{\Theta_n}) \Vert \\
&\leq & \Vert g_{\Theta} - \rho_n^{-1}(g_{\Theta}) \Vert + \Vert
\rho_n^{-1}\Vert \Vert (\chi_{\Theta}- \chi_{\Theta_n})g \Vert
\end{eqnarray*}
and that implies the statement by $(w_n,v_n) \to (1,0)$ and strong
continuity of the representation. (ii) If, on the other hand,
$\kappa_- = \kappa_+$, i.e. the sequence $\Theta_n$ tends to a
single point, we have, again by contractivity of $\mu^2
R(\Delta_{\Theta_n},\mu^2)$,
\begin{equation*}
    \Vert\mu^2 R(\Delta_{\Theta_n},\mu^2)g_{\Theta_n}\Vert = \Vert\mu^2
    R(\Delta_{\Theta_n},\mu^2)g_{\Theta_n}\Vert\leq\Vert g_{\Theta_n}\Vert
\end{equation*}
which tends to zero as $n$ tends to infinity.
\end{Proof}

Since the sets $\PP_n (U)$ are discrete, we do not have to prove a corresponding result for fixed $n\neq\infty$. It
remains to show convergence for partition solvers to partitions $p_n\in\PP_n(U)$ that converge to a partition $p\in\PP
(U)$ which is not contained in $\PP_n (U)$ for any $n$.

\begin{Lemma}\label{psolverdrei} Let $p_n\in\PP_n(U)$ be a
sequence of partitions converging to $p\in\PP(U)$ with respect to Hausdorff metric. Then
\begin{equation*}
    \lim_{n\to\infty} f_{\mu,n,g}^* (p_n) = f_{\mu,\infty,g}^*
    (p).
\end{equation*}
\end{Lemma}

\begin{Proof} Recall the notations from Corollary \ref{RPMS},\ref{RPBZ},
respectively. By Lemma \ref{convchar}, the proof consists of
considering the following two cases: (i) Let
$\Theta=(\kappa_-,\kappa_+)\in\iota (p)$ and
$\Theta_n\in\iota(p_n)$ such that
$\overline{\Theta_n}\to\overline{\Theta}$. Then $\Theta_n =
({\scriptstyle \frac{\kappa_{r_n}}{n}} ,
{\scriptstyle\frac{\kappa_{r_n+1}}{n}})$ and we have to show that
for $g\in L^2 (U)$
\begin{equation*}
    -\mu^2 \tau_n \,
    R(n^2 B(r_n),\mu^2)\,\pi_n g^t\to -\mu^2
    R(\Delta_{\Theta},\mu^2)g
\end{equation*}
as $n$ tends to infinity. We observe that, by cancelling $-\mu^2$, this can be done by proving {\em strong resolvent
convergence} of the operator
\begin{equation*}
    \Delta_n := \tau_n\, n^2 B(r_n)\circ\pi_n
\end{equation*}
to the Neumann laplacian $\Delta_{\Theta}$. That follows from
\begin{eqnarray*}
    \mu^{-2}(\Delta_n - \mu^2) &=& \mu^{-2}(\tau_n \,n^2 B (r_n)\, \pi_n - \mu^2 \tau_n\,\pi_n +
    \mu^2 (\tau_n\,\pi_n - 1 ) )\\
    &=& \tau_n \, \mu^{-2}(n^2 B (r_n)-\mu^2)\, \pi_n  -
     (1 - \tau_n\,\pi_n )
\end{eqnarray*}
where $\tau_n \,\pi_n$ is nothing else but conditional expectation
with respect to $\sigma_n$ on the subinterval $\lbrack
{\scriptstyle \frac{\kappa_{r_n}}{n}},
{\scriptstyle\frac{\kappa_{r_{n+1}-1}}{n}} )$. Hence for all
$f\in\DD (\Delta_n)$, we have by $\pi_n (1 - \tau_n\,\pi_n)=0$
\begin{eqnarray*}
    & & \tau_n\,
    R(n^2 B(r_n),\mu^2)\,\pi_n (\Delta_n - \mu^2)f \\
    &=& \tau_n\,\pi_n
    f - \mu^2 \tau_n\,
    R(n^2 B(r_n),\mu^2)\,\pi_n (1 - \tau_n\,\pi_n)f \\
    &=& \tau_n\,\pi_n f.
\end{eqnarray*}
Together with $\tau_n\,\pi_n f \to f_{\Theta}$ this implies that
the resolvents $R(\Delta_n,\mu^2)$ tend to the same limit as the
operators $\tau_n\, R(n^2 B(r_n),\mu^2)\,\pi_n$ strongly as $n$
tends to infinity.

It remains to prove strong resolvent convergence of $\Delta_n$ to
$\Delta_{\Theta}$. By \cite{Kato:80}, Corollary 1.6, p. 429, and
since $\lbrace f\in C^{\infty}(\overline{\Theta}):
f^{\prime}\vert_{\partial \overline{\Theta}} = 0 \rbrace$ forms a
core for the Neumann laplacian, we just have to prove that
\begin{equation*}
    \lim_{n\to\infty} \Vert\Delta_n f - \Delta_{\Theta} f \Vert = 0
\end{equation*}
for all $f$ in the core. Let thus
\begin{equation*}
    f(x) = f (x_0) + (x-x_0) f^{\prime} (x_0) +
    \frac{1}{2}f^{\prime\prime}(x_0) (x-x_0)^2 + R(x)
\end{equation*}
be a Taylor expansion of $f$ with remainder $\vert R(x)\vert =
O(\vert x\vert^3)$ uniformly for $x$ in a compact interval
containing $\overline{\Theta}$. Then, controlling the error by
Taylor expansion yields
\begin{eqnarray*}
    \Delta_n f
    &=& \sum_{\kappa = \kappa_{r_n}+1}^{\kappa_{r_n+1}-2}
    (f^{\prime\prime}(\kn)+ O(1/n))\chi_{\lbrack \kn,\kpn)} \\
    & & - (n f^{\prime}(\kappa_-) + O(1))\chi_{\lbrack{\scriptstyle \frac{\kappa_{r_n}}{n}},{\scriptstyle \frac{\kappa_{r_n}+1}{n}}
    )}
    + (n f^{\prime}(\kappa_+) + O(1))\chi_{\lbrack{\scriptstyle \frac{\kappa_{r_{n+1}}-1}{n}},{\scriptstyle \frac{\kappa_{r_{n+1}}}{n}}
    )}.
\end{eqnarray*}
Hence, the boundary conditions
$f^{\prime}(\kappa_+)=f^{\prime}(\kappa_-)=0$ imply strong
convergence of $R(\Delta_n,\mu^2)$ to $R(\Delta_{\Theta},\mu^2)$
uniformly on compact sets $\mu \in \lbrack a,b \rbrack\subset
\R^+$, $b>a>0$, contained in the {\em resolvent set}. By Lipschitz
continuity established in Lemma \ref{psolvereins} this extends to
uniform convergence of $-\mu^2 R(\Delta_n,\mu^2)g$ to $-\mu^2
R(\Delta_{\Theta},\mu^2)g$ for $\mu \in \lbrack 0, b\rbrack$. (2)
For a sequence of intervals collapsing to a point, uniform
contractivity $\Vert \mu^2 R(\Delta_n,\mu^2)\Vert\leq 1$ implies
convergence of the corresponding partition solvers to $0\in L^2
(U)$ as in the proof of the preceding lemma.
\end{Proof}

\section{Compactness of the Set of Minimizers}\label{CP}

By the reduction principle, the set of partition solvers to a
fixed function $g\in L^2(U)$ contains the set of minimizers of
$F$. Now we will prove, that the set of partition solvers is
compact, which implies the same for the set of minimizers. Let
$g\in L^2 (U)$, $b>0$ and
\begin{equation*}
    \MM_b (g) := \lbrace f^{*}_{\mu,n,g}(p)\, :\, p\in\PP_n (U),\mu\in \lbrack 0,b\rbrack,n\in\N\cup\lbrace\infty\rbrace\rbrace .
\end{equation*}
In this section, will prove
\begin{Theorem}\label{mainzwei} For all $b>0$ and $g\in L^2 (U)$,
the set $\MM_b (g)\subset  L^2 (U)$ is compact.
\end{Theorem}

This result implies the existence of global minimizers.

\begin{Theorem}[Existence of Minimizers]\label{mini} Let $\gamma ,\mu\geq 0$ and $u\in L^2(U)$. Then:
\begin{enumerate}
\item The set of global minimizers $f^*$ of $\MS_{\gamma,\mu,g}$
is non-empty. \item The set of global minimizers $f_n^*$ of
$\BZ_{\gamma,\mu,g}^n$ is non-empty.
\end{enumerate}
\end{Theorem}

\begin{Proof} (i) Recall that the reduced Mumford - Shah functional is given by
\begin{equation*}
    \ms_{\gamma,\mu,g}(p) = \gamma j(p) + \Vert g \Vert^2 - \langle g,f^*_{\mu,\infty,g}(p)\rangle .
\end{equation*}
By Theorem \ref{mainzwei}, minimizing $\MS_{\gamma,\mu,g}$ on $L^2
(U)$ is hence equivalent to the minimization of
\begin{equation*}
    \gamma j(p(f)) + \Vert g \Vert^2 - \langle f,g\rangle
\end{equation*}
for functions $f$ in the compact set $\MM_{\mu}(g)$. For $\gamma =
0$, this function is continuous and hence assumes its minimum on
the compact set. For $\gamma > 0$, the function is lower semi -
continuous by Corollary \ref{cardinality} and assumes hence as
well its minimum on a compact set. (ii) The set $\PP_n (U)$ is
finite, hence the minimization can be reduced to the minimization
with respect to a finite number of partition solvers which implies
that the functional assumes its minimum.
\end{Proof}

Now we start with the proof of Theorem \ref{mainzwei}.

\subsection{A Preliminary Lemma}

First of all, we will prove a lemma that will simplify the
discussion considerably. It states that a proof of compactness for
all $g\in L^2(U)$ can be reduced to a proof for a {\em total
subset} of signals in $L^2 (U)$. Denote therefore by
\begin{equation*}
    \MM_b := \lbrace g\in L^2 (U)\,:\,\MM_b (g)\subset L^2
    (U)\,\,\mathrm{\, compact}\rbrace
\end{equation*}
the set of all signals, such that $\MM_b (g)$ is compact in $L^2
(U)$. Then we have the following statement:

\begin{Lemma}\label{bigsimplification} For all $b>0$, the set $\MM_b\subset L^2
(U)$ is a closed linear subspace.
\end{Lemma}

\begin{Proof} (1) Let $g_1,g_2\in \MM_b$, $\lambda_1,\lambda_2\in\R$. Then, by linearity of
the partition solver
\begin{equation*}
f^{*}_{\mu,\lambda_1 g_1+\lambda_2 g_2,n}(p)=\lambda_1
f^{*}_{\mu,g_1,n}(p)+ \lambda_2 f^{*}_{\mu,g_2,n}(p),
\end{equation*}
and hence $\MM_b (\lambda_1 g_1 + \lambda_2 g_2)\subset
\lambda_1\MM_b (g_1) + \lambda_2\MM_b (g_2)$. The right hand side
is compact by continuity of the linear operations on $L^2 (U)$.
(2) Let $g_k\in\MM_b$ be a sequence of signals converging to $g\in
L^2 (U)$. The mappings $g\mapsto f^{*}_{\mu,n,g}(p)$ are
contractions. Hence
\begin{equation*}
\Vert f^{*}_{\mu,n,g_k}(p) - f^{*}_{\mu,n,g}(p)\Vert\leq \Vert g_k
- g\Vert,
\end{equation*}
i.e. the sets $\MM_b (g_k)$ converge to $\MM_b (g)$ as closed
subsets of $L^2(U)$ with respect to Hausdorff metric. That implies
that $\MM_b (g)$ is compact, and therefore $\MM_b$ is closed.
\end{Proof}
The preceding lemma implies immediately, that if $\MM_b$ contains
a {\em total set}, i.e. a set of functions such that its closed
linear hull equals $L^2 (U)$, it already equals the whole space.
Such a set is provided by the {\em Heaviside functions} on $U$.

\subsection{Compactness}

Denote by
\begin{equation*}
    \HH (U) := \lbrace \chi_a := \chi_{\lbrack a,1\rbrack}\,:\,
    0\leq a < 1 \rbrace \subset L^2 (U)
\end{equation*}
the set of {\em Heaviside functions} on the interval $U$. The
linear hull of $\HH (U)$ is provided by the step functions which
are dense in $L^2 (U)$, hence $\HH (U)$ forms a total set. The key
result for compactness reads

\begin{Proposition} For all $\chi_a\in \HH (U)$ we have
\begin{equation*}
    \MM (\chi_a) \subset \lbrace f\in\MM_b (\chi_a)\,:\, \vert
    p(f)\vert\leq 6\rbrace ,
\end{equation*}
and $\MM_b (\chi_a)$ is compact.
\end{Proposition}

\begin{Proof} Let $p =\lbrace 0 < x_1 < ...< x_r <1\rbrace$ be an
arbitrary partition such that $a\in \lbrack x_k,x_{k+1})$. Then,
$\chi_a$ and $\chi_{a,n}:= E(\chi_a \vert\sigma_n)$ are constant
on $\lbrack 0, x_k)$ (equal to zero) and $\lbrack x_{k+1}, 1)$
(equal to one). Since the kernels of $\Delta_{\Theta},\Delta_n$
are provided by constant functions, vectors, respectively, the
respective resolvents map the conditional expectation $\chi_{a,n}$
onto functions which are constant (and equal to zero) on $\lbrack
0, x_{n-1})$ and constant (equal to one) on $\lbrack
x_{k+2},1\rbrack$. That implies, they have at most six jumps
points, contained in the set $\lbrace
x_{k-1},x_k,x_{k+1},x_{k+2}\rbrace\cup\lbrace 0,1\rbrace$. By
Lemma \ref{toppartitions} (2) and Corollary \ref{cardinality}, the
set $\lbrace p\in\PP (U)\,:\, \vert p \vert\leq 6\rbrace$ is
compact as a closed subset of a compact set and hence the closed
subset
\begin{equation*}
    \C := \E (\T) \cap \T\times \lbrace p\in\PP (U)\,:\, \vert p \vert
    \leq 6\rbrace
\end{equation*}
is as well compact. That implies finally by Theorem \ref{maineins}
compactness of
\begin{equation*}
\MM_b (\chi_a) \subset f^{*}(\lbrack 0, b\rbrack\times \C).
\end{equation*}
\end{Proof}

\begin{Remark} Recall that by Definition \ref{discdef}, the notion
of jump differs from the notion of discontinuity in the discrete
setting.
\end{Remark}

By Lemma \ref{bigsimplification}, that completes the proof of
Theorem \ref{mainzwei}.

\section{$\mathbf{\Gamma}$-Convergence and Weak Derivatives}\label{GC}
The notion of $\Gamma$-convergence is important for the
investigation of minimizers. It provides some rather general
sufficient condition for the possibility to approximate minimizers
of a given functional by minimizers of a sequence of approximating
functionals.

\begin{Definition} Let $X$ be a metric space and $F_j:X\to \overline{\R}$
be a sequence of functions. Then, $F_j$ converges to $F$ in the
sense of $\Gamma$-convergence -- in the sequel we will write
shortly $F_j \stackrel{\Gamma}{\rightarrow} F_{\infty}$ -- if
\begin{enumerate}
\item For all $x\in X$ and all sequences $x_j \stackrel{X}{\to}x$
we have
\begin{equation}\label{gammaeins}
F_{\infty}(x)\leq \liminf_{j\to\infty} F_j (x_j).
\end{equation}
\item For all $x\in X$ there is a sequence $\hat{x}_j
\stackrel{X}{\to}x$ such that
\begin{equation}\label{gammazwei}
F_{\infty}(x)\geq \limsup_{j\to\infty} F_j (\hat{x}_j).
\end{equation}
\end{enumerate}
\end{Definition}
Essentially, $\Gamma$-convergence is important due to the
following facts (cf. Theorem 5.3.6 of \cite{A:Bee:93})

\begin{Theorem}\label{MainGamma} Suppose $F_j \stackrel{\Gamma}{\rightarrow} F_{\infty}$ and denote by $\mathrm{argmin} F$ the set of minimizers of
$F$.Then
\begin{enumerate}
  \item  For any converging sequence $x_{j,j\in\N}$, $x_j\in\argmin \,F_j$, we have necessarily $\lim_{j\to\infty}x_j\in\argmin \,F_\infty$.
  \item If there is a compact subset $K\subset X$ such that $\emptyset\neq\argmin \,F_j\subset K$ for large enough $j$, then $\argmin \,F_\infty\ne\emptyset$ and
    \begin{equation*}
      \lim_{j\to\infty} d(x_j,\argmin \,F_\infty) = 0
    \end{equation*}
for any  sequence $x_{j,j\in\N}$, $x_j\in\argmin F_j$. \item If,
additionally, $\argmin F_\infty$ is a singleton $\lbrace
x_\infty\rbrace$ then
\begin{equation*}
      \lim_{j\to\infty} x_j=x_\infty
\end{equation*}
for any  sequence $x_{j,j\in\N}$ with $x_j\in\argmin F_j$.
\end{enumerate}
\end{Theorem}
In the sequel, we will prove that the Blake - Zisserman
functionals $\BZ^n$ converge to the associated Mumford - Shah
functionals $\MS$ in $\Gamma-$sense  as $n$ tends to infinity. By
the compactness results established in the preceding section, this
implies convergence of the associated minimizers.

The crucial step will be the understanding of the behavior of the
Blake - Zisserman penalty $\Phi^n_{\gamma,\mu}$ (see Definition
\ref{BZD}) as $n$ tends to infinity. For reasons that will become
clear in the sequel, we prove a parameter dependent result.

\begin{Proposition}\label{penale} Let $\gamma_n$, $\mu_n$
sequences of non - negative numbers converging to $\gamma, \mu
\geq 0$. Then the following two statements are valid which imply
$\Gamma$- convergence of $\Phi^n_{\gamma_n,\mu_n,g}$ to the
Mumford - Shah penalty to parameters $\gamma,\mu$:
\begin{enumerate}
\item Let $f=t + F\in \SBV_2 (U)$. Then there is a sequence
$\tilde{f}_{n,n\in\N}$ converging to $f$ in $L^2 (U)$ such that
\begin{equation*}
    \limsup_{n\to\infty} \sum_{\kappa =0}^{n-2} \min\lbrace \koenn
(\tilde{f}_n (\kpn) - \tilde{f}_n (\kn))^2,\gamma_n\rbrace \leq
\gamma j(f) + \mu^{-2} \int_0^1 \vert f^{\prime}(x)\vert^2 dx.
\end{equation*}
\item Let $f\in L^2 (U)$. For all sequences $f_{n,n\in\N}$
converging to $f$ in $L^2 (U)$ we have
\begin{eqnarray*}
& &\liminf_{n\to\infty} \sum_{\kappa =0}^{n-2} \min\lbrace \koenn
(f_n (\kpn) - f_n (\kn))^2,\gamma_n\rbrace \\ & &
\geq\left\lbrace\begin{array}{ll} \gamma j(t) + \mu^{-2} \int_0^1
\vert
f^{\prime}(x)\vert^2 dx & f\in\SBV_2 (U)\\
\infty & \mathrm{else}\end{array}\right. .
\end{eqnarray*}
\end{enumerate}
Note that in the case $\mu = 0$, the right hand side of both
inequalities is only finite for $f\in\TT (U)\subset \SBV_2 (U)$,
since $f^{\prime}\equiv 0$ exactly for those functions.
\end{Proposition}

To prove Proposition \ref{penale}, we have to collect several
facts about the relation of weak differentiation and approximation
by step functions. This will be done in the following two
subsection. The final proof of Proposition \ref{penale} is given
in \ref{PPen}.

As an introductory step, we discuss weak differentiability in $L^2
(U)$.

\subsection{A Characterization of Weak Differentiability}

A function $f$ is Sobolev - differentiable on the one-sphere, i.e.
$f\in H^1 (S^1)$, if and only if $f\in L^2 (S^1)$ and its {\em
Fourier coefficients} $\hat{f}(k),k\in\Z$ with respect to the
orthonormal base $\phi_k (t) := e^{2\pi i kt}$ satisfy
\begin{equation}\label{fouriersphere}
    \sum_{k\in\Z} k^2\,\vert\hat{f}(k)\vert^2 < \infty
\end{equation}
(see for instance \cite{Roe:88}, Section 3). Although $L^2 (U)$
and $L^2 (S^1)$ can be identified by considering $S^1$ as the
identification $U/\lbrace 0,1\rbrace$, the situation is different
for the Sobolev spaces $H^1(\Omega)$ and $H^1 (S^1)$ ($\Omega =
(01,)$, cf. Definition \ref{SBV}). By {\em Sobolev's embedding
theorem}, $H^1 (\Omega)\subset C(U)$, and $f\in H^1 (\Omega)$ can
be identified with some $\tilde{f}\in H^1 (S^1)$ if and only if
$f(0) = f(1)$ for the continuous representative of $f$. Hence,
$H^1 (S^1)\subset H^1 (\Omega)$ is a linear subspace of
codimension one. In order to prove $\Gamma$-convergence for the
Blake-Zisserman functional, the first step will be to find a
suitable analogue of the characterization (\ref{fouriersphere})
for $H^1 (\Omega)$.

\begin{Lemma}\label{charsobolew} Equivalent are
\begin{enumerate}
\item $f\in H^1 (\Omega)$, \item there is some $\alpha\in\R$ such
that $f-g_{\alpha}\in H^1 (S^1)$ where $g_{\alpha} (t) = \alpha
t$, \item there is some $\alpha\in\R$ such that
\begin{equation*}
    \sum_{k\in\Z} \left\vert 2\pi i k \hat{f}(k) -
    \alpha\right\vert^2 < \infty .
\end{equation*}
where $\hat{f}(k)$ denote the Fourier coefficients of $f$,
\end{enumerate}
In case $f\in H^1 (\Omega)$ we have for $\alpha\in\R$ as above
\begin{equation*}
    \sum_{k\in\Z} \left\vert 2\pi i k \hat{f}(k) -
    \alpha\right\vert^2 = \int_0^1 \vert (f - g_{\alpha})^{\prime}\vert^2
    dt = \int_0^1 \vert f^{\prime}\vert^2 dt - (f(1) - f(0))^2,
\end{equation*}
hence $\alpha = f(1) - f(0)$.
\end{Lemma}

\begin{Proof} (1) First, we show that (i) implies (ii). Let $f\in H^1 (\Omega)$. By Sobolev's embedding
theorem, there is a representative $f_0$ of $f$ such that $f_0\in
C(U)$. Let $\phi\in C^{\infty}(U)$ with $d^k\phi/dt^k (0) =
d^k\phi/dt^k (1)$ for all $k=0,1,...$ and $\alpha := f(1) - f(0)$.
Then
\begin{equation*}
    \int_0^1 (f^{\prime}-\alpha ) \phi dt = \phi (0) (f-g_{\alpha})\vert_0^1
    - \int_0^1 (f - g_{\alpha}) \phi^{\prime} dt
\end{equation*}
and hence $f^{\prime}-\alpha$ is the weak derivative of
$f-g_{\alpha}$ in $H^1 (S^1)$. Since $U$ is compact,
square-integrability of $f^{\prime}$ implies square-integrability
of $f^{\prime}-\alpha$. (2) To see that (ii) implies (i), let
$f-g_{\alpha}\in H^1 (S^1)$ and $\phi\in
C^{\infty}_0(\Omega)\subset C^{\infty}(S^1)$. Then, since
$g_{\alpha}\in H^1 (\Omega)$,
\begin{equation*}
   \int_0^1 f\phi^{\prime} dt=\int_0^1 (f-g_{\alpha} + g_{\alpha} )\phi^{\prime} dt
   =  -\int_0^1 ((f-g_{\alpha})^{\prime} + \alpha
   )\phi dt .
\end{equation*}
Hence, $f$ is weakly differentiable on $\Omega$ with square
integrable weak derivative $f^{\prime}=(f-g_{\alpha})^{\prime} +
\alpha$. (3) The equivalence of (ii) and (iii) follows from
\begin{equation*}
    \widehat{f^{\prime}} (k) = \int_0^1 dt e^{-2\pi ikt}
    f^{\prime}(t) = f(1)-f(0) + 2\pi i k\hat{f}(k)
\end{equation*}
and the equivalence of (i) and (ii).
\end{Proof}

\begin{Lemma}\label{mollifier} Let a family of {\em mollifiers} be given by
\begin{equation*}
    h_a(x) := \left\lbrace\begin{array}{ll} \frac{1}{a}x & 0<x\leq a \\ 1 & a<x\leq 1-a \\
    \frac{1}{a}(1-x) & 1-a\leq x < 1\end{array}\right.
\end{equation*}
where $0<a<1/2$. Then, if $f\in L^2 (U)$, we have $fh_a\in L^2
(U)$ for all $a>0$ with
\begin{equation*}
    \Vert fh_a\Vert\leq \Vert f\Vert
\end{equation*}
and $fh_a$ converges to $f$ in $L^2 (U)$ as $a$ tends to zero.
\end{Lemma}

\begin{Proof} Clearly $0\leq h_a \leq 1$ and $h_a \to 1$ in $L^2
(U)$ as $a$ tends to zero. Hence
\begin{equation*}
    \Vert fh_a\Vert^2= \int_0^1 \vert f\vert^2 h_a^2 dt\leq \int_0^1 \vert f\vert^2  dt=\Vert
    f\Vert^2
\end{equation*}
and
\begin{equation*}
    \Vert f(h_a-1)\Vert^2\leq \int_0^a \vert f\vert^2  dt + \int_{1-a}^1 \vert f\vert^2  dt
\end{equation*}
which tends to zero as $a$ tends to zero.
\end{Proof}
These simple observations are the starting point for the
constructions in the next sections.

\subsection{Approximation and Weak Differentiability}

Given $f^{(n)}\in \TT_n (U)$, we consider
\begin{equation*}
F^{(n)} (x) := n (f^{(n)} (x+1/n) - f^{(n)}(x))
\end{equation*}
where we understand $f^{(n)}$ to be continued periodically to the
real axis. Then
\begin{equation*}
\int_0^{1-\frac{1}{n}} F^{(n)}(s)^2 ds = \sum_{\kappa =0}^{n-2} n
(f^{(n)}(\kpn) - f^{(n)}(\kn))^2 .
\end{equation*}
Now we rewrite the integral in terms of the Fourier coefficients
of $F^{(n)}$ with respect to the orthonormal base
\begin{equation*}
    \phi_k^{(n)} (t) := \left\lbrace\begin{array}{ll}\sqrt{\frac{n}{n-1}}\,e^{\frac{2\pi i
    nkt}{n-1}}& ,0< t \leq 1-\frac{1}{n} \\ 0 & ,1-\frac{1}{n} < t <
    1 \end{array}\right.
\end{equation*}
on $L^2 (\lbrack 0, 1-\frac{1}{n}\rbrack)$. As $n$ tends to
infinity, the individual base vectors tend to $\phi_k (t) :=
e^{2\pi i kt}$ and even the Fourier coefficients converge as we
will show in the next lemma.

\begin{Lemma}\label{coeffasymp} Let $f^{(n)}_n\in\N$ be a sequence of step functions with $f^{(n)}\in\TT_n (U)$ converging to $f$
in $L^2 (U)$, and $F^{(n)}$ as above. For fixed $k\in \Z$, $a\in
(0,1/2)$ we have
\begin{equation}
    \lim_{n\to\infty} \langle F^{(n)}h_a, \phi^{(n)}_k -
    \phi_k\rangle = 0 .
\end{equation}
\end{Lemma}

\begin{Proof} (1) First of all, we show that
\begin{equation}
    \lim_{n\to\infty} \langle F^{(n)}h_a, \phi^{(n)}_k -
    \tilde{\phi}^{(n)}_k\rangle = 0
\end{equation}
where
\begin{equation*}
\tilde{\phi}^{(n)}_k (t) := \sqrt{\frac{n}{n-1}}\,e^{\frac{2\pi i
    nkt}{n-1}}.
\end{equation*}
Assume without loss of generality that $n$ is so large that $1/n <
a$. Then by Cauchy - Schwarz inequality
\begin{eqnarray*}
\vert\langle F^{(n)}h_a, \phi^{(n)}_k -
\tilde{\phi}^{(n)}_k\rangle\vert &=&
\vert\langle F^{(n)},h_a \chi_{\lbrack 1-\frac{1}{n},1)}
\tilde{\phi}^{(n)}_k\rangle\vert \leq \Vert F^{(n)}\Vert \,\Vert h_a \chi_{\lbrack 1 - \frac{1}{n},1)} \tilde{\phi}^{(n)}_k\Vert \\
&\leq& 2n \Vert f^{(n)}\Vert\,\sqrt{\frac{n\int_{1-\frac{1}{n}}^1
dt (1-t)^2}{a^2(n-1)}}\leq \frac{2\,\Vert
f^{(n)}\Vert}{\sqrt{3(n-1)}} .
\end{eqnarray*}
This expression tends to zero since, by assumption, $\Vert
f^{(n)}\Vert\to\Vert f \Vert <\infty$. (2) Now we have by
substitution of $u:= t + 1/n$, using periodicity and $\vert h_a
\vert \leq 1$
\begin{eqnarray*}
    & & \vert\langle F^{(n)}h_a, \tilde{\phi}^{(n)}_k -
    \phi_k\rangle\vert \\
    &\leq& n\int_0^{1} dt \vert(f^{(n)} (t+1/n) - f^{(n)}(t))h_a(\tilde{\phi}^{(n)}_k  - \phi_k)(t)\vert \\
    &\leq& n\int_0^{1} du \vert f^{(n)} (u) (\tilde{\phi}^{(n)}_k (u-1/n) - \tilde{\phi}^{(n)}_k(u)
    + \phi_k (u) - \phi_k(u-1/n))\vert.
\end{eqnarray*}
But now
\begin{eqnarray*}
& & \vert\tilde{\phi}^{(n)}_k (u-1/n) - \tilde{\phi}^{(n)}_k(u)+
\phi_k (u) - \phi_k(u-1/n)\vert^2 \\
&=& \left\vert\sqrt{\frac{n}{n-1}}(1-e^{\frac{2\pi i
kt}{n-1}})(1-e^{-\frac{2\pi i kt}{n}})e^{2\pi i kt}\right\vert^2
\\
&=& \frac{4n}{n-1}\left(1-\cos\frac{2\pi
kt}{n-1}\right)\left(1-\cos\frac{2\pi i kt}{n}\right) =
O({\scriptstyle\frac{1}{n^4}}).
\end{eqnarray*}
Hence
\begin{equation*}
    \vert\langle F^{(n)}h_a, \tilde{\phi}^{(n)}_k -
    \phi_k\rangle\vert = O({\scriptstyle\frac{1}{n}})
\end{equation*}
tends to zero as $n$ tends to infinity.
\end{Proof}

\begin{Lemma}\label{fourcoeff} Let $f^{(n)}_n\in\N$ be a sequence of step functions with
$f^{(n)}\in\TT_n (U)$ converging to $f$ in $L^2 (U)$, and
$F^{(n)}$ as above. Let $\hat{f}(k)$ be the Fourier expansion of
$f$ with respect to the orthonormal base $\phi_k (t) := e^{2\pi i
kt}$ and $h_a$ as above. Let
\begin{equation}\label{defect}
    \alpha_k (a,n) := \frac{1}{a}\left(\int_{1-a}^1 f^{(n)}\phi_k(t)dt - \int_0^{a} f^{(n)}\phi_k(t)dt\right)
\end{equation}
and $\alpha_k(a):=\lim_{n\to\infty} \alpha_k (a,n)$. Then
\begin{equation*}
    \lim_{n\to\infty}\widehat{F^{(n)}h_a}(k) =   2\pi i k\widehat{f\,h_a}(k) +
     \alpha_k (a).
\end{equation*}
\end{Lemma}

\begin{Proof} The function $h_a\phi_k$ is everywhere {\em left-differentiable} with left differential
\begin{equation*}
    (h_a\phi_k)^{\prime}_{-}(t) := \lim_{s\to 0} \frac{h_a\phi_k (t) - h_a\phi_k (t-s)}{s}.
\end{equation*}
Hence, by substituting $t+1/n=u$, we have
\begin{eqnarray*}
    \widehat{F^{(n)}h_a}(k) &=& n\int_0^1 (f^{(n)}(t+1/n)-f^{(n)}(t))h_a\phi_k(t)dt \\
    &=& -n\int_0^1 f^{(n)}(u)(h_a\phi_k(u) - h_a\phi_k (u-1/n))du
    .
\end{eqnarray*}
The function $h_a \phi_k$ is uniformly Lipschitz with Lipschity
constant $C_{a,k}$. Hence $\vert n(h_a\phi_k(u) - h_a\phi_k
(u-1/n))\vert\leq C_{a,k}$ and thus by dominated convergence
\begin{eqnarray*}
    & &\lim_{n\to \infty}\widehat{F^{(n)}h_a}(k)\\
    &=& -\lim_{n\to \infty}\int_0^1 f^{(n)}(h_a\phi_k)_{-}^{\prime}du \\
    &=& -\lim_{n\to \infty}\int_0^1 f^{(n)}(2\pi i k h_a - h_{a-}^{\prime})\phi_k\,du \\
    &=& - 2\pi i k\lim_{n\to \infty}\int_0^1 f^{(n)} h_a \phi_k du + \lim_{n\to \infty}\int_0^1 f^{(n)}\frac{1}{a}
    (\chi_{\lbrack 0,a)} - \chi_{\lbrack 1-a,1)})\phi_k du .
\end{eqnarray*}
Convergence $f^{(n)}\to f$ finally implies the statement.
\end{Proof}

\begin{Lemma}\label{nocheinsunddanngutenacht} Let $\alpha_k (a)$
be as above. Let $a_m\to 0$ be a sequence of positive numbers such
that $\lim_{m\to\infty} \alpha_k (a_m) = \alpha_k \neq \pm\infty$
exists for some $k\in\Z$. Then this limit exists for all $k\in\Z$
and equals $\alpha_0$.
\end{Lemma}

\begin{Proof} $f\in L^2 (U)$ implies $f\in L^1 (U)$. On the other
hand, there is a constant $C_k>0$ such that $\vert 1 - \phi_k
(t)\vert\leq C_k \min \lbrace\vert t \vert,\vert 1-t \vert\rbrace$
for all $t\in \lbrack 0,a\rbrack\cup \lbrack 1-a,1\rbrack$. Using
(\ref{defect}), the statement follows now from
\begin{eqnarray*}
    \vert \alpha_0(a) - \alpha_k(a)\vert&\leq& \frac{1}{a}\left\lbrack\int_0^a
    \vert f (1-\phi_k)\vert dt + \int_{1-a}^1 \vert f
    (1-\phi_k)\vert  dt\right\rbrack \\ &\leq& C_k\,\left\lbrack\int_0^a
    \vert f\vert dt + \int_{1-a}^1 \vert f
    \vert  dt\right\rbrack
\end{eqnarray*}
which tends to zero as $a\to 0$.
\end{Proof}
From this considerations, we obtain, having in mind inequality
(\ref{gammaeins}):

\begin{Lemma}\label{liminf} Let $f^{(n)}\in\TT_n (U)$ a sequence of step functions converging to $f$
in $L^2 (U)$, $F^{(n)}$ as above. Then we have
\begin{equation*}
    \liminf_{n\to\infty}\sum_{\kappa =0}^{n-2} n (f^{(n)}(\kpn) - f^{(n)}(\kn))^2 \geq
    \left\lbrace\begin{array}{ll}\int_0^1 \vert f^{\prime}(x)\vert^2 dx
    & \mathrm{if}\, f\in H^1(\Omega) \\ \infty &
    \mathrm{if}\, f\notin H^1 (\Omega)\end{array}\right. .
\end{equation*}
\end{Lemma}

\begin{Proof} We have by Lemma \ref{mollifier} for all $a>0$
\begin{equation*}
 \int_0^{1-\frac{1}{n}} (F^{(n)}(x))^2 dx\geq \int_0^{1-\frac{1}{n}} (F^{(n)}h_a(x))^2 dx = \sum_{k\in\Z}
\left\vert\langle F^{(n)}h_a,\phi_k^{(n)}\rangle \right\vert^2
\end{equation*}
and all summands are non-negative. Hence by {\em Fatou's Lemma}
and Lemma \ref{coeffasymp} and \ref{fourcoeff}, we have for all
$a>0$:
\begin{eqnarray*}
    & & \liminf_{n\to\infty}\sum_{k\in\Z}
\left\vert\langle F^{(n)}h_a,\phi_k^{(n)}\rangle \right\vert^2 \geq \sum_{k\in\Z}\liminf_{n\to\infty}
    \left\vert\langle F^{(n)}h_a,\phi_k^{(n)}\rangle \right\vert^2 \\
    &=& \sum_{k\in\Z}\liminf_{n\to\infty}
    \left\vert\widehat{F^{(n)}h_a}(k)\right\vert^2 = \vert \alpha_0 (a)\vert^2 + \sum_{k\neq 0}
    \left\vert 2\pi i k\widehat{fh_a}(k)-
    \alpha_k (a)\right\vert^2 .
\end{eqnarray*}
By Lemma \ref{nocheinsunddanngutenacht}, either $\liminf_{a\to
0}\vert \alpha_0 (a)\vert=\infty$, or for all sequences $a_m\to 0$
for which the limit $\lim_{m\to \infty}\vert \alpha_0 (a_m)\vert^2
= \vert a_0\vert^2$ exists, we have as well $\lim_{m\to
\infty}\vert \alpha_k (a_m)\vert^2 = \vert a_0\vert^2$ for all
$k\in\Z$. Taking such a subsequence, we obtain by {\em Fatou's
Lemma} and Lemma \ref{mollifier}
\begin{eqnarray*}
    & & \liminf_{n\to\infty} \int_0^{1-\frac{1}{n}} (F^{(n)}(x))^2
    dx \\
    &\geq& \liminf_{m\to \infty}\left\lbrack \vert \alpha_0 (a_m)\vert^2 + \sum_{k\neq 0}
    \left\vert 2\pi i k\widehat{fh_{a_m}}(k)-
    \alpha_k (a_m)\right\vert^2 \right\rbrack \\
    &\geq&  \lim_{m\to \infty}\vert \alpha_0 (a_m)\vert^2 + \sum_{k\neq 0}
    \lim_{m\to\infty}\left\vert 2\pi i k\widehat{fh_{a_m}}(k)-
    \alpha_k (a_m)\right\vert^2 \\
    &\geq&  \lim_{m\to \infty}\vert \alpha_0 (a_m)\vert^2 + \sum_{k\neq 0}
    \lim_{m\to\infty}\left\vert 2\pi i k\widehat{f}(k)-
    \alpha_k (a_m)\right\vert^2 \\
    &=&  \vert \alpha_0 \vert^2 + \sum_{k\neq 0}
    \left\vert 2\pi i k\widehat{f}(k)-
    \alpha_0\right\vert^2 .
\end{eqnarray*}
By Lemma \ref{charsobolew}, $f\in H^1(\Omega)$ if and only if
there is some $\alpha\in\R$ such that the sum on the right hand
side is finite. Hence, if $f\notin H^1(\Omega)$, the right hand
side is always infinite. If $f\in H^1 (\Omega)$, $f$ has a
continuous version and thus $\lim_{a\to 0}\alpha_0 (a) = f(0) -
f(1)$ which implies by Lemma \ref{charsobolew}, that the limes
inferior equals $\int \vert f^{\prime} \vert^2 dt$.
\end{Proof}

\noindent The result corresponding to inequality (\ref{gammazwei})
reads as follows:

\begin{Lemma}\label{limsup} Let $f\in H^1(\Omega)$, $\widetilde{f^{(n)}} :=
E(f\vert\sigma_n)$ the conditional expectation with respect to the
sigma-algebra $\sigma_n$ and $\tilde{F}^{(n)}=
n(\widetilde{f^{(n)}}(x+1/n)-\widetilde{f^{(n)}}(x))$as above.
Then
\begin{equation*}
    \limsup_{n\to\infty} \int_0^{1-\frac{1}{n}} \vert \tilde{F}^{(n)}(x)\vert^2 dx
    \leq\int_0^1 \vert f^{\prime}(x)\vert^2  dx .
\end{equation*}
\end{Lemma}

\begin{Proof} By the definition of conditional expectation and {\em Jensen's inequality} (applied
to the probability measure $n\,dx$ on $\lbrack \kmn,\kn)$) we obtain
\begin{eqnarray*}
& &\sum_{\kappa =0}^{n-2} n (\widetilde{f^{(n)}}(\kn) -
\widetilde{f^{(n)}}(\kmn))^2 = \sum_{\kappa =0}^{n-2} n
\left\lbrack n
\int_{\kmn}^{\kn} (f(x+1/n) - f(x))dx\right\rbrack^2 \\
&\leq& \sum_{\kappa =0}^{n-2}  \int_{\kmn}^{\kn}(n
 (f(x+1/n) - f(x))^2dx = \int_{0}^{1-\frac{1}{n}}(n (f(x+1/n) - f(x))^2 dx.
\end{eqnarray*}
By {\em Lebesgue's differentiation theorem}
(\cite{Wheeden:Zygmund:1977}, (7.2) Theorem, p. 100) we have
convergence
\begin{equation*}
    n (f(x+1/n) - f(x))= \frac{1}{\lambda(\lbrack x,x+1/n))}\int_{x}^{x+1/n} f^{\prime}(u)
    du
\end{equation*} to $f^{\prime}(x)$
for Lebesgue-almost all $x\in I$. On the other hand
\begin{equation*}
\vert n (f(x+1/n) - f(x))\vert\leq  2 S^{*}(x):= 2
\sup_{V}\frac{1}{\vert V\vert}\int_{V} \vert f^{\prime}(u)\vert
du,
\end{equation*}
where $V$ is any open sub-interval $V\subset I$ such that $x\in V$. The function $S^{*}$ is not integrable, except for
$f=0$ almost surely (see \cite{Wheeden:Zygmund:1977}, p. 105). But by the Lemma of {\em Hardy-Littlewood}
(\cite{Wheeden:Zygmund:1977}, (7.9) Theorem, p. 105) there is a constant $c>0$, such that
\begin{equation*}
\lambda(\lbrace x\in U\,:\, S^{*} (x) > \alpha\rbrace)\leq
\frac{c}{\alpha}\int_{U}\vert f^{\prime}(u)\vert du < \infty,
\end{equation*}
since $L^2(U)\subset L^1 (U)$. That implies
\begin{eqnarray*}
    \lambda (\lbrace x\in U\,:\, \vert n (f(x+1/n) - f(x))\vert >
    \alpha\rbrace )
    &\leq&\lambda( \lbrace x\in U\,:\, 2 S^{*} (x) >
    \alpha\rbrace ) \\
    &\leq&\frac{2 c}{\alpha}\int_{U}\vert f^{\prime}(u)\vert du
\end{eqnarray*}
and the sequence of difference functions $n (f(x+1/n) - f(x))$ is
therefore {\em uniformly integrable}. Thus, we may interchange
limit and integration which implies the statement.
\end{Proof}

\subsection{Convergence of the Smoothness Penalty}\label{PPen}

As a consequence of the considerations in the preceding
subsection, we prove now Proposition \ref{penale}, the
corresponding $\Gamma$-convergence result for the Blake -
Zisserman penalty $\Phi^n_{\gamma,\mu,u}$. Note that the case $\mu
= 0$ requires some care.

\begin{Proof} (i) Let $f=F+t\in \SBV_2 (U)$. From the sequence $f_{n,n\in\N}$, $f_n\in\TT_n(U)$ we construct the
decomposition $f^F_n := f_n - f_n^t$ where
\begin{equation*}
    f_n^t := \sum_{\kappa =0}^{n-1}t(\kn) \chi_{\lbrack\kn,\kpn)}\in\TT_n(U).
\end{equation*}
We thus have $f_n = f_n^F + f_n^t$ and $f_n^t\to t$ by boundedness
of $t$, hence as well $f_n^F\to F$. Furthermore $f_n^t (\kpn) -
f_n^t (\kn) = 0$ if the interval $( \kn,\kpn\rbrack$ contains no
jump of the step function $t$ and there are only finitely many
intervals that contain a jump, namely at most $j(t)$. Hence
\begin{eqnarray*}
& & \sum_{\kappa =0}^{n-2} \min\lbrace \koenn (f_n (\kpn) - f_n
(\kn))^2,\gamma_n\rbrace \\
&& = \sum_{\lbrace\kappa : ( \kn ,\kpn \rbrack\cap
p(t)=\emptyset\rbrace} \min\lbrace \koenn (f_n
(\kpn ) - f_n (\kn ))^2,\gamma_n\rbrace \\
&&+ \sum_{\lbrace k: ( \kn, \kpn\rbrack\cap
p(t)\neq\emptyset\rbrace} \min\lbrace \koenn (f_n (\kpn) -
f_n (\kn))^2,\gamma_n\rbrace\\
&& = \sum_{\lbrace \kappa : ( \kn ,\kpn\rbrack\cap
p(t)=\emptyset\rbrace} \min\lbrace \koenn (f_n^F
(\kpn) - f_n^F (\kn))^2,\gamma_n\rbrace \\
&&+ \sum_{\lbrace \kappa : ( \kn ,\kpnzwei\rbrack\cap
p(t)\neq\emptyset\rbrace} \min\lbrace \koenn
(f_n^F (\kpn) - f_n^F (\kn)+ f_n^t (\kpn) - f_n^t (\kn))^2,\gamma_n\rbrace\\
\end{eqnarray*}
Now we consider the sequence of conditional expectations
$\tilde{f}_n^F := E(f - t\,\vert \,\sigma_n)$. Then $\tilde{f}_n^F
:= E(F\,\vert\,\sigma_n)$ and we have by Lemma \ref{limsup}
\begin{eqnarray*}
&& \limsup_{n\to\infty}\sum_{\lbrace \kappa : ( \kn
,\kpn\rbrack\cap p(t)=\emptyset\rbrace} \min\lbrace \koenn
(\tilde{f}_n^F
(\kpn) - \tilde{f}_n^F (\kn))^2,\gamma_n\rbrace \\
&\leq& \lim_{n\to\infty}\sum_{\kappa =0}^{n-2} \koenn
(\tilde{f}_n^F (\kpn) - \tilde{f}_n^F (\kn))^2 = \mu^{-2}_n
\int_0^1 \vert f^{\prime}(x)\vert^2 dx .
\end{eqnarray*}
For the analysis of the exceptional intervals, we use the fact
that $F$ is {\em absolutely continuous} (see
\cite{Wheeden:Zygmund:1977}, p. 115). Therefore, for all $\epsilon
> 0$ there is an $n_0$, such that for all $n\geq n_0$ and $\kappa =0,...,n-1$ we have
\begin{equation*}
    \sup_{x,x^{\prime}\in \lbrack\kn , \kpn)}\vert
F(x)-F(x^{\prime})\vert\leq \epsilon.
\end{equation*}
By the contraction property of conditional expectation that
implies for all $n\geq n_0$, $\kappa =0,...,n-1$ that
$\vert\tilde{f}_n^F (\kpn) - \tilde{f}_n^F (\kn)\vert \leq 2
\epsilon$. Choose now $n_0$ so large that $\epsilon < \delta/4$
where $\delta := \min_{x\in J(t)} \vert t(x)-t(x^-)\vert$ is the
heigt of the smallest jump of the step function and additionally,
such that all exceptional intervals contain exactly one
discontinuity of $t$. That implies for the exceptional intervals
\begin{equation*}
    \vert \tilde{f}_n^F (\kpn) - \tilde{f}_n^F (\kn)+ f_n^t (\kpn) - f_n^t (\kn)\vert \geq \delta/4.
\end{equation*}
Thus for all $n$ with $n\geq n_0$ and $n\delta^2/16\mu^2_n > \gamma_n$ we have
\begin{eqnarray*}
&& \sum_{\lbrace \kappa : ( \kn ,\kpn \rbrack\cap
p(t)\neq\emptyset\rbrace} \min\lbrace \koenn
(f_n^F (\kpn) - f_n^F (\kn)+ f_n^t (\kpn) - f_n^t (\kn))^2,\gamma_n\rbrace \\
&=& \sum_{\lbrace \kappa : ( \kn ,\kpn \rbrack\cap
p(t)\neq\emptyset\rbrace} \min\lbrace n \delta^2/16\mu^2_n
,\gamma_n\rbrace = \gamma_n\,j(t).
\end{eqnarray*}

\noindent (ii) Let $f\in L^2(U)$, $f_n\in \TT_n(U)$ with $f_n\to
f$. Then
\begin{eqnarray*}
& & \sum_{\kappa =0}^{n-2} \min\lbrace \koenn (f_n (\kpn) - f_n
(\kn))^2,\gamma_n\rbrace \\&=& \sum_{k\in a_n^c}
 \koenn (f_n (\kpn) - f_n (\kn))^2 + \sum_{\kappa\in a_n} \gamma_n \\
&=& \gamma_n\,\vert a_n\vert + \sum_{\kappa\in a_n^c}
 \koenn (f_n (\kpn) - f_n (\kn))^2 ,
\end{eqnarray*}
where $a_n := \lbrace \kappa\leq n-2: \min\lbrace \koeff (g_n
(\kpn) - g_n (\kn))^2,\gamma\rbrace=\gamma\rbrace$ and
$a_n^{\mathrm c} := \lbrace 0,...,n-2\rbrace - a_n$. Without loss
of generality, we consider the subsequence $f_n^{\mathrm{inf}}$ of
$f_n$ with
\begin{eqnarray*}
    & &\liminf_{n\to\infty}\sum_{\kappa =0}^{n-2} \min\lbrace \koenn (f_n (\kpn) - f_n
(\kn))^2,\gamma_n\rbrace \\&=& \lim_{n\to\infty}\sum_{\kappa
=0}^{n-2} \min\lbrace \koenn (f_n^{\mathrm{inf}} (\kpn) -
f_n^{\mathrm{inf}} (\kn))^2,\gamma_n\rbrace
\end{eqnarray*}
with corresponding exceptional sets $a_n^{\mathrm{inf}}$. The sets
$p_n^{\mathrm{inf}} := \lbrace \kappa /n : \kappa\in
a_n^{\mathrm{inf}}\rbrace \cup \lbrace 0,1\rbrace$ are finite and
hence closed in $\lbrack 0,1\rbrack$. Passing to another
subsequence of $f_n^{\mathrm{inf}}$ if necessary, the compactness
of the set of closed subsets of $\lbrack 0,1\rbrack$ with respect
to Hausdorff distance implies, that the sequence of sets
$p_n^{\mathrm{inf}}$ converges to a closed subset $c\subset
\lbrack 0,1\rbrack$. By (i), convergence to a function $f\in\SBV_2
(U)$ implies by the absolute continuity of $f$ off the jump points
that the exceptional set contains only finitely many points. Thus
$\vert c\vert=\infty$ implies $f\notin \SBV_2 (U)$.

\noindent Therefore assume $\vert c \vert = K < \infty$. Hence
$c\in\PP (U)$. In that case, Lemma \ref{liminf} yields the
following alternative: Either $f\in \SBV_2 (U)$, then $f\in
H^1(\Theta)$ for all $\Theta\in\iota (c)$ and the limit of the
subsequence is greater or equal to
\begin{equation*}
    \gamma K + \mu^{-2} \sum_{\Theta \in\iota (c)} \vert f^{\prime}_{\Theta}\vert^2
    dx
\end{equation*}
where $f^{\prime}_{\Theta} = f^{\prime}\vert_{\Theta}$, or
$f\notin\SBV_2 (U)$ which implies that the limit is $\infty$. That
implies the statement, in particular for $\mu = 0$.
\end{Proof}

\section{Dependence on the Parameters}

In the final section, we will prove the Theorem \ref{Hauptsatz}.
According to Theorem \ref{MainGamma}, the proof follows from
$\Gamma$-convergence together with the fact -- already established
in Section \ref{CP} -- that the set of minimizers is compact.
Thus, we start by showing that the functionals in question depend
continuously on the respective parameters in an appropriate sense.

\subsection{$\mathrm{\Gamma}$-Continuity of the Segmentation Family}

We consider the three-dimensional {\em (pseudo-) cube} $\QQ$ given
by
\begin{equation*}
    \QQ := \R^+_0\times \R^+_0 \times \T
\end{equation*}
and the corresponding family of functionals $F(q) : L^2 (U)\to \R$, $q\in\QQ$ defined in Theorem \ref{Hauptsatz}. The
statement about $\Gamma$-continuity on the cube now reads as follows:

\begin{Theorem}\label{GammaCube} Let $q_{s,s\in\N}$ with $q_s := (\gamma_s,\mu_s,t_s)\in \QQ$ a sequence of
parameters converging to $q:=(\gamma,\mu,t)\in \QQ$. Then
\begin{equation*}
    F (q_s)\stackrel{\Gamma}{\rightarrow} F (q)
\end{equation*}
as $s$ tends to infinity.
\end{Theorem}

To prove this, the crucial point is the statement about
$\Gamma$-convergence of the penalizers established in Proposition
\ref{penale}. However, we will need two additional lemmas, the
first of which states that the discrete $L^2$-distance used in the
Blake - Zisserman functional converges to the continuous
$L^2$-distance.

\begin{Lemma}\label{distance} Let $f,g\in L^2 (U)$, $f_{n,n\in\N}$ with $f_n \in
\TT_n(U)$ be a sequence of step functions converging to $f$ in
$L^2 (U)$ and $g_{n,n\in\N}$ with $g_n := E(g\,\vert\,\sigma_n)$
the sequence of conditional expectations. Then
\begin{equation*}
    \lim_{n\to\infty}\frac{1}{n}\sum_{\kappa =0}^{n-1} (f_n -
g_n)^2(\kn) = \Vert f - g \Vert^2
\end{equation*}
\end{Lemma}

\begin{Proof} First of all, $E(f_n \,\vert\, \sigma_n )=f_n$ by a standard property of
conditional expectation. Furthermore, {\em martingale convergence}
(see e.g. \cite{Williams:91}, Ch. 12) implies
\begin{equation*}
    \lim_{n\to\infty} E(f_n - f\,\vert \, \sigma_n) = 0.
\end{equation*}
Hence
\begin{eqnarray*}
&&\lim_{n\to\infty}\frac{1}{n}\sum_{\kappa =0}^{n-1} (f_n -
g_n)^2(\kn) = \lim_{n\to\infty}\int_0^1 E(f_n - g \vert\sigma_n)^2
dx \\ &=& \lim_{n\to\infty}\int_0^1 (E(f_n - f \vert\sigma_n) +
E(f-g\,\vert \sigma_n))^2 dx  = \lim_{n\to\infty}\int_0^1
E(f-g\,\vert \sigma_n)^2 dx  \\
\end{eqnarray*}
Again by {\em martingale convergence}, $E(f-g\,\vert \sigma_n)$
tends to $f - g$.
\end{Proof}

By the second lemma, we prove a result about approximation of
$L^2$-functions under additional constraints about the location of
the jumps of the approximating step functions.

\begin{Lemma}\label{dasletzte} Let $f\in L^2 (U)$ and $k_m$ a
sequence of non negative integers such that $k_m\leq m$ and
$k_m\to\infty$. Then there is a sequence of functions
\begin{equation*}
    f_{m}\in \lbrace f\in\TT_m (U)\,:\, j(f_{m})\leq
    k_m\rbrace
\end{equation*}
such that $f_{m}\to f$ in $L^2 (U)$.
\end{Lemma}

\begin{Proof} We consider three families associated to $f$:
\begin{enumerate}
\item $F_k:=\argmin_{t\in\TT(U),j(t)\leq k}\Vert t-f\Vert$,
$k=1,2,...$ -- the minimum is not necessarily unique. \item
$F_{k,m}:= E(F_k\,\vert\,\sigma_m)$ -- we have $j(F_{k,m})\leq
2k$, $k,m=1,2,...$. \item $f_{k,m}:=\argmin_{t\in\TT_m(U),j(t)\leq
k}\Vert t-f\Vert$, $k,m=1,2,...$ -- again, the minimum is not
necessarily unique.
\end{enumerate}
By $k_m\to\infty$, we have for all $K>0$ some $M_0> 0$ such that
$k_m>2K$ for all $m\geq M_0$. Then
\begin{equation*}
    \Vert f_{k_m,m} - f\Vert \leq \Vert F_{K,m} - f \Vert \leq
    \Vert F_{K,m} - F_K\Vert + \Vert F_{K} - f\Vert .
\end{equation*}
By construction
\begin{equation*}
    \lim_{m\to\infty}\Vert F_{K,m} - F_K\Vert =0, \lim_{K\to\infty}\Vert F_{K} - f\Vert = 0
\end{equation*}
and thus
\begin{equation*}
    \limsup_{m\to\infty}\Vert f_{k_m,m} - f\Vert \leq\Vert F_{K} - f\Vert
\end{equation*}
arbitrarily small and non negative. That implies the statement.
\end{Proof}

Now we can prove the theorem stated above.

\begin{Proof} For the proof, we have to consider three different cases depending
on the location of $q$.

\noindent{\sc 1st case.} $\gamma,\mu\geq 0,t > 0$: In that case $q_s\to q$ if and only if $\gamma_s\to\gamma,
\mu_s\to\mu$ and $t_s = t:= 1/n$ for $s\geq s_0$. That implies for $s\geq s_0$
\begin{eqnarray*}
    & & F (\gamma_s,\mu_s,t_s)(f) - F (\gamma,\mu,t)(f) \\
    &=& \BZ_{\gamma_s,\mu_s,u}^n (f) - \BZ_{\gamma,\mu,u}^n
    (f) \\
    &=&\left\lbrace\begin{array}{ll}\sum_{\kappa =0}^{n-2}
    \min\lbrace \koess(f(\kpn) - f(\kn))^2 , \gamma_s
    \rbrace & f\in \TT_n(U)\\- \min\lbrace \koeff(f(\kpn) - f(\kn))^2 , \gamma
    \rbrace& \\ 0 & \mathrm{else} \end{array}\right. .
\end{eqnarray*}
For step functions $f\in\TT_n(U)$, $\Vert f\Vert_{L^2 (U)} \leq C$ implies $\sup \vert f\vert \leq C\sqrt{n}$. Thus $F
(\gamma_s,\mu_s,t_s)$ converges to $F (\gamma,\mu,t)$ uniformly on balls $\lbrace \Vert f \Vert_{L^2(U)} \leq
C\rbrace$. That implies $\Gamma$-convergence of the functionals.

\noindent{\sc 2nd case.} $\gamma>0,\mu\geq 0,t = 0$: In that case, $\Gamma$-convergence follows from Proposition
\ref{penale} for the penalizer and Lemma \ref{distance} for the distance term.

\noindent{\sc 3rd case.} $\gamma=0,\mu\geq 0, t=0$: In that case
we have $n_s:=1/t_s\to\infty$ and for $f_s\in\TT_{n_s}(U)$
\begin{equation*}
    \hat{d}^{n_s}_g (f_s)\leq \BZ_{\gamma_s,\mu_s,g}^{n_s} (f_s)
\leq n_s\gamma_s + \hat{d}^{n_s}_g (f_{s}).
\end{equation*}
Let now $f\in L^2(U)$ and $f_s\in\TT_{n_s}(U)$ a sequence of
functions with $f_s\to f$ in $L^2(U)$. Then
\begin{equation*}
    \liminf_{s\to\infty} \BZ_{\gamma_s,\mu_s,g}^{n_s} (f_s)\geq\liminf_{s\to\infty}
\hat{d}^{n_s}_g (f_s) = \lim_{s\to\infty} \hat{d}^{n_s}_g (f_s)=
d_g (f).
\end{equation*}
That is condition (\ref{gammaeins}) for $\Gamma-$convergence. Let
now $\hat{f}_s\in\TT_{n_s}(U)$ a sequence of step functions such
that each $\hat{f}_s$ is some best approximation of $f$ in the set
$M_s:= \lbrace t\in\TT_n(U)\,:\, j(t)\leq 1/\gamma_s \rbrace$.
Then $\gamma_s j(\hat{f}_s) \to 0$, $\hat{f}_s\to f$ and by Lemma
\ref{dasletzte}, this sequence fulfills condition
(\ref{gammazwei}).
\end{Proof}

\subsection{Proof of the Main Theorem}

The proof of Theorem \ref{Hauptsatz} is now finally a consequence
of Theorem \ref{MainGamma} together with the compactness of the
set of minimizers: Existence (i) of the minimizers is provided by
Theorem \ref{mini}. Convergence (iia) of the minimizers follows
from Theorem \ref{MainGamma}, (i) together with Theorem
\ref{GammaCube}, the result on $\Gamma$-convergence established
above. For all sequences $q_s = (\gamma_s,\mu_s,t_s)$ converging
to $q=(\gamma,\mu,t)\in\QQ$, there is some $b > 0$ with $\mu_s\leq
b$ for all $s$. Since, by the reduction principle, all possible
minimizers are thus contained in the set $\MM_b (g)$ (see Section
\ref{CP}), the existence of a convergent subsequence (iib) follows
from the compactness result Theorem \ref{mainzwei}.

\bibliographystyle{plain}
\bibliography{MumfordBib}

\end{document}